\documentclass{article}
\usepackage{amsmath,amsthm,amsopn,amstext,amscd,amsfonts,amssymb}
\usepackage{natbib}

\def \tr {\mathop{\rm tr}\nolimits}
\def \re {\mathop{\rm Re}\nolimits}

\def \Vol {\mathop{\rm Vol}\nolimits}

\def \etr {\mathop{\rm etr}\nolimits}
\def \diag {\mathop{\rm diag}\nolimits}
\def \build#1#2#3{\mathrel{\mathop{#1}\limits^{#2}_{#3}}}

\renewenvironment{abstract}
                 {\vspace{6pt}
                  \begin{center}
                  \begin{minipage}{5in}
                  \centerline{\textbf{Abstract}}
                  \noindent\ignorespaces
                 }
                 {\end{minipage}\end{center}}
\newtheorem{thm}{\textbf{Theorem}}[section]
\newtheorem{cor}{\textbf{Corollary}}[section]
\newtheorem{lem}{\textbf{Lemma}}[section]
\theoremstyle{definition}

\newtheorem{rem}{\textbf{Remark}}[section]

\title{\Large \textbf{Random matrix theory and multivariate statistics}}
\author{
  \textbf{Jos\'e A. D\'{\i}az-Garc\'{\i}a} \thanks{Corresponding author\newline
    }\\
  {\normalsize Department of Statistics and Computation} \\
  {\normalsize 25350 Buenavista, Saltillo, Coahuila, Mexico} \\
  {\normalsize E-mail: jadiaz@uaaan.mx} \\[2ex]
  \textbf{Ram\'on Guti\'errez J\'aimez} \\
  {\normalsize Department of Statistics and O.R.} \\
  {\normalsize University of Granada} \\
  {\normalsize Granada 18071, Spain}\\
  {\normalsize E-mail: rgjaimez@ugr.es}\\
}
\date{}
\begin{document}
\maketitle

\begin{abstract}
Some tools and ideas are interchanged between random matrix theory and multivariate
statistics. In the context of the random matrix theory, classes of spherical and
generalised Wishart random matrix ensemble, containing as particular cases the
classical random matrix ensembles, are proposed. Some properties of these classes of
ensemble are analysed.  In addition, the random matrix ensemble approach is extended
and a unified theory proposed for the study of distributions for real normed division
algebras in the context of multivariate statistics.
\end{abstract}

\section{Introduction}\label{sec1}

Random $\,$ matrices $\,$ first $\,$appeared in studies of mathematical statistics in
the 1930s but did not attract much attention at the time. However, in time the
foundations were laid for one of the main areas of present-day statistics,
\emph{Multivariate Statistics} or \emph{Multivariate Statistical Analysis}. The
advances made have been compiled in excellent books such as \citet{ro:57} and
\citet{a:58}, later in \citet{sk:79} and \citet{m:82} and recently in \citet{gn:00}
and \citet{ea:07}, among many other important texts. In general, these writings are
based on the multivariate normal (Gaussian) distribution.  Nevertheless, over the
last 20 years, a new current has been developed and consolidated  to constitute what
is now known as \emph{Genaralised Multivariate Analysis}, in which the multivariate
normal distribution is replaced by a family of elliptical distributions, containing
as particular cases the Normal, Student-$t$, Pearson type VII, Logistic and Kotz
distributions, among many others, see \citet{fkn:90}, \citet{fz:90} and
\citet{gv:93}. This family contains many distributions with heavier or lighter tails
than the Gaussian distribution, thus offering a more flexible framework for modelling
phenomena and experiments in more general situations.

Random matrices are not only of interest to statisticians; they are studied in many
disciplines of science, engineering and finance. This approach has been consolidated
and is generally known as \emph{Random Matrix Theory}, see \citet{me:91},
\citet{f:05}, \citet{rva:05a} and \citet{er:05}, among many others. Fundamentally,
this theory is concerned with the following question: consider a large matrix whose
elements are random variables with given probability laws. Then, what can be said
about the probabilities of a few of its eigenvalues or of a few of its eigenvectors?
This question is relevant to our understanding of the statistical behaviour of slow
neutron resonances in nuclear physics, a field in which the issue was first addressed
in the 1950s and where it is intensively studied by physicists. The question was
later found to be important in other areas of physics and mathematics, such as the
characterisation of chaotic systems, the elastodynamic properties of structural
materials, conductivity in disordered metals, the distribution of the values of the
Riemann zeta function on the critical line, the enumeration of permutations with
certain particularities, the counting of certain knots and links, quantum gravity,
quantum chromo dynamics, string theory, and others, see \citet{me:91} and
\citet{f:05}. This approach has also motivated the study of random linear systems
(operators), random matrix calculus (calculus of Jacobians of matrix factorisations)
and numerical stochastic algorithms, among other mathematical fields, see
\citet{er:05}. A very interesting characteristic of the results obtained in random
matrix theory is that they are vaild for real, complex and quaternionic random
matrices.

From the standpoint of the multivariate statistics, the main focus of statistical
studies has been developed for real random matrices. However, many texts have
proposed the extension of some of these tests to the case of complex random matrices,
see for example \citet{w:56}, \citet{g:63}, \citet{j:64}, \citet{k:65},
\citet[Section 6.6, p. 295]{ha:70}, \citet{m:97} and \citet{mdm:06}, among many
others. Studies examining the case of quaternionic random matrices case are much less
common, see \citet{b:00} and \citet{lx:09}.

Form the point of view of the statistical and random matrix theory, the relevance of
\emph{the octonions} is not clear. An excellent review of the history, construction
and many other properties of the octonions is \citet{b:02}, where it is stated that:
\begin{center}
\begin{minipage}[t]{4.5in}
\begin{sl}
``Their relevance to geometry was quite obscure until 1925, when \'Elie Cartan
described `triality' -- the symmetry between vector and spinors in 8-dimensional
Euclidian space. Their potential relevance to physics was noticed in a 1934 paper by
Jordan, von Neumann and Wigner on the foundations of quantum mechanics...Work along
these lines continued quite slowly until the 1980s, when it was realized that the
octionions explain some curious features of string theory... \textbf{However, there
is still no \emph{proof} that the octonions are useful for understanding the real
world}. We can only hope that eventually this question will be settled one way or
another."
\end{sl}
\end{minipage}
\end{center}

Others problems, such as the existence of real or imaginary eigenvalues for $2 \times
2$ and $3 \times 3$ octonionic Hermitian matrices are treated in \citet{dm:98} and
\citet{dm:99}.  For the sake of completeness, in the present study the case of the
octonions is considered, but the veracity of the results in this case can only be
conjectured. In particular, \citet[Section 1.4.5, pp. 22-24]{f:05} it is proveed that
the bi-dimensional density function of the eigenvalue, for a Gaussian ensemble of a
$2 \times 2$ octonionic matrix, is obtained from the general joint density function
of the eigenvalues for the Gaussian ensemble, assuming $m = 2$ and $\beta = 8$, see
Section \ref{sec3}.

In the present study, we propose an interchange between some ideas and tools of
random matrix theory and multivariate statistics. Specifically, in Section \ref{sec2}
general expressions are given for the Jacobians of certain matrix transformation,
together with the factorisation for real normed division algebras, a question that is
fundamental in the study of general random matrix variate distributions.
Vector-spherical, spherical, vector-spherical-Laguerre and spherical-Laguerre random
matrix ensembles are studied in Section \ref{sec3}, where it is observed that the
classical random matrix ensembles are obtained as particular cases of the above
general random matrix variate distributions. Finally, the density function of the
matrix variate normal, Wishart, elliptical, generalised Wishart and beta type I and
II distributions are studied for real normed division algebras.

\section{ Jacobians}\label{sec2}

In this section we discuss the formulas containing the Jacobians of familiar matrix
transformations and factorisations using the parameter $\beta$ (defined below).
First, however, some notation must be established.

The parameter $\beta$ has traditionally been used to count the real dimension of the
underlying normed division algebra. In other branches of mathematics, the parameter
$\alpha = 2/\beta$ is used.

\medskip

\begin{table}[!h]
  \centering
  \caption{Values of $\beta = 2/\alpha$ parameter.}\label{table1}
\begin{tabular}{l|l l}
  \hline\hline
  $\beta$ & $\alpha$ & Normed divison algebra \\
  \hline
  1 & 2 & real ($\Re$) \\
  2 & 1 & complex ($\mathfrak{C}$) \\
  4 & 1/2 & quaternionic ($\mathfrak{H}$) \\
  8 & 1/4 & octonion ($\mathfrak{O}$) \\
  \hline
\end{tabular}
\end{table}

\begin{rem}
In \citet{b:02} the following remarks are made:
\begin{center}
\begin{minipage}[t]{4.5in}
\begin{sl}
``There are exactly four normed division algebras: the real numbers ($\Re$), complex
numbers ($\mathfrak{C}$), quaternions ($\mathfrak{H}$), and octonions
($\mathfrak{O}$). The real numbers are the dependable breadwinner of the family, the
complete ordered field we all rely on. The complex numbers are a slightly flashier
but still respectable younger brother: not ordered, but algebraically complete. The
quaternions, being noncommutative, are the eccentric cousin who is shunned at
important family gatherings. But the octonions are the crazy old uncle nobody lets
out of the attic: they are nonassociative."
\end{sl}
\end{minipage}
\end{center}
In addition, take into account that $\Re$, $\mathfrak{C}$, $\mathfrak{H}$ and
$\mathfrak{O}$ are the only normed division algebras; moreover, they are the only
alternative division algebras, and all division algebras have a real dimension of $1,
2, 4$ or $8$, see \citet[Theorems 1, 2 and 3]{b:02}.
\end{rem}

Let ${\mathcal L}^{\beta}_{m,n}$ be the linear space of all $n \times m$ matrices of
rank $m \leq n$ over $\mathfrak{F}$ with $m$ distinct positive singular values, where
$\mathfrak{F}$ denotes a \emph{real finite-dimensional normed division algebra}. Let
$\mathfrak{F}^{n \times m}$ be the set of all $n \times m$ matrices over
$\mathfrak{F}$. The dimension of $\mathfrak{F}^{n \times m}$ over $\Re$ is $\beta
mn$. Let $\mathbf{A} \in \mathfrak{F}^{n \times m}$, then $\mathbf{A}^{*} =
\overline{\mathbf{A}}^{T}$ denotes the usual conjugate transpose. The set of matrices
$\mathbf{H}_{1} \in \mathfrak{F}^{n \times m}$ such that
$\mathbf{H}_{1}^{*}\mathbf{H}_{1} = \mathbf{I}_{m}$ is a manifold denoted ${\mathcal
V}_{m,n}^{\beta}$, termed \emph{Stiefel manifold} ($\mathbf{H}_{1}$ is also known as
\emph{semi-orthogonal} ($\beta = 1$), \emph{semi-unitary} ($\beta = 2$),
\emph{semi-symplectic} ($\beta = 4$) and \emph{semi-exceptional type} ($\beta = 8$)
matrices, see \citet{dm:99}). The dimension of $\mathcal{V}_{m,n}^{\beta}$ over $\Re$
is $[\beta mn - m(m-1)\beta/2 -m]$. In particular, ${\mathcal V}_{m,m}^{\beta}$ with
dimension over $\Re$, $[m(m+1)\beta/2 - m]$, is the maximal compact subgroup
$\mathfrak{U}^{\beta}(m)$ of ${\mathcal L}^{\beta}_{m,m}$ and consist of all matrices
$\mathbf{H} \in \mathfrak{F}^{m \times m}$ such that $\mathbf{H}^{*}\mathbf{H} =
\mathbf{I}_{m}$. Therefore, $\mathfrak{U}^{\beta}(m)$ is the \emph{real orthogonal
group} $\mathcal{O}(m)$ ($\beta = 1$), the \emph{unitary group} $\mathcal{U}(m)$
($\beta = 2$), \emph{compact symplectic group} $\mathcal{S}p(m)$ ($\beta = 4$) or
\emph{exceptional type matrices} $\mathcal{O}o(m)$ ($\beta = 8$), for $\mathfrak{F} =
\Re$, $\mathfrak{C}$, $\mathfrak{H}$ or $\mathfrak{O}$, respectively. Denote by
${\mathfrak S}_{m}^{\beta}$ the real vector space of all $\mathbf{S} \in
\mathfrak{F}^{m \times m}$ such that $\mathbf{S} = \mathbf{S}^{*}$. Let
$\mathfrak{P}_{m}^{\beta}$ be the \emph{cone of positive definite matrices}
$\mathbf{S} \in \mathfrak{F}^{m \times m}$; then $\mathfrak{P}_{m}^{\beta}$ is an
open subset of ${\mathfrak S}_{m}^{\beta}$. Over $\Re$, ${\mathfrak S}_{m}^{\beta}$
consist of \emph{symmetric} matrices; over $\mathfrak{C}$, \emph{Hermitian} matrices;
over $\mathfrak{H}$, \emph{quaternionic Hermitian} matrices (also termed
\emph{self-dual matrices}) and over $\mathfrak{O}$, \emph{octonionic Hermitian}
matrices. Generically, the elements of $\mathfrak{S}_{m}^{\beta}$ are termed as
\textbf{Hermitian matrices}, irrespective of the nature of $\mathfrak{F}$. The
dimension of $\mathfrak{S}_{m}^{\beta}$ over $\Re$ is $[m(m-1)\beta+2]/2$. Let
$\mathfrak{D}_{m}^{\beta}$ be the \emph{diagonal subgroup} of
$\mathcal{L}_{m,m}^{\beta}$ consisting of all $\mathbf{D} \in \mathfrak{F}^{m \times
m}$, $\mathbf{D} = \diag(d_{1}, \dots,d_{m})$. Let $\mathfrak{T}_{U}^{\beta}(m)$ be
the subgroup of all \emph{upper triangular} matrices $\mathbf{T} \in \mathfrak{F}^{m
\times m}$ such that $t_{ij} = 0$ for $1 < i < j \leq m$; and let
$\mathfrak{T}_{L}^{\beta}(m)$ be the opposed \emph{lower triangular} subgroup
$\mathfrak{T}_{L}^{\beta}(m) = \left(\mathfrak{T}_{U}^{\beta}(m)\right)^{T}$. For any
matrix $\mathbf{X} \in \mathfrak{F}^{n \times m}$, $d\mathbf{X}$ denotes the\emph{
matrix of differentials} $(dx_{ij})$. Finally, we define the \emph{measure} or volume
element $(d\mathbf{X})$ when $\mathbf{X} \in \mathfrak{F}^{m \times n},
\mathfrak{S}_{m}^{\beta}$, $\mathfrak{D}_{m}^{\beta}$ or $\mathcal{V}_{m,n}^{\beta}$,
see \citet{d:02}.

If $\mathbf{X} \in \mathfrak{F}^{n \times m}$ then $(d\mathbf{X})$ (the Lebesgue
measure in $\mathfrak{F}^{n \times m}$) denotes the exterior product of the $\beta
mn$ functionally independent variables
$$
  (d\mathbf{X}) = \bigwedge_{i = 1}^{n}\bigwedge_{j = 1}^{m}\bigwedge_{k =
  1}^{\beta}dx_{ij}^{(k)}.
$$
\begin{rem}
\begin{itemize} Note that for $x_{ij} \in \mathfrak{F}$
$$
    dx_{ij} = \bigwedge_{k = 1}^{\beta}dx_{ij}^{(k)}.
  $$
In particular for $\mathfrak{F} = \Re$, $\mathfrak{C}$, $\mathfrak{H}$ or
$\mathfrak{O}$ we have
  \item $x_{ij} \in \Re$ then
  $$
    dx_{ij} = \bigwedge_{k = 1}^{1}dx_{ij}^{(k)}= dx_{ij}.
  $$
  \item $x_{ij} = x_{ij}^{(1)} + ix_{ij}^{(2)} \in \mathfrak{C}$, then
  $$
    dx_{ij} = dx_{ij}^{(1)} \wedge dx_{ij}^{(2)} = \bigwedge_{k = 1}^{2}dx_{ij}^{(k)}.
  $$
  \item $x_{ij} = x_{ij}^{(1)} + ix_{ij}^{(2)}+ jx_{ij}^{(3)}+ kx_{ij}^{(4)} \in
  \mathfrak{H}$, then
  $$
    dx_{ij} = dx_{ij}^{(1)} \wedge dx_{ij}^{(2)}\wedge dx_{ij}^{(3)}\wedge dx_{ij}^{(4)}
    = \bigwedge_{k = 1}^{4}dx_{ij}^{(k)}.
  $$
  \item $x_{ij} = x_{ij}^{(1)} + e_{1}x_{ij}^{(2)}+ e_{2}x_{ij}^{(3)}+ e_{3}x_{ij}^{(4)} +
  e_{4}x_{ij}^{(5)} + e_{5}x_{ij}^{(6)}+ e_{6}x_{ij}^{(7)}+ e_{7}x_{ij}^{(8)} \in
  \mathfrak{O}$, then
  $$
    dx_{ij} = dx_{ij}^{(1)} \wedge dx_{ij}^{(2)} \wedge dx_{ij}^{(3)} \wedge dx_{ij}^{(4)}
    \wedge dx_{ij}^{(5)} \wedge dx_{ij}^{(6)} \wedge dx_{ij}^{(7)} \wedge
    dx_{ij}^{(8)} = \bigwedge_{k = 1}^{8}dx_{ij}^{(k)}
  $$
\end{itemize}
\end{rem}
If $\mathbf{S} \in \mathfrak{S}_{m}^{\beta}$ (or $\mathbf{S} \in
\mathfrak{T}_{L}^{\beta}(m)$) then $(d\mathbf{S})$ (the Lebesgue measure in
$\mathfrak{S}_{m}^{\beta}$ or in $\mathfrak{T}_{L}^{\beta}(m)$) denotes the exterior
product of the $m(m+1)\beta/2$ functionally independent variables (or denotes the
exterior product of the $m(m-1)\beta/2 + n$ functionally independent variables, if
$s_{ii} \in \Re$ for all $i = 1, \dots, m$)
$$
  (d\mathbf{S}) = \left\{
                    \begin{array}{ll}
                      \displaystyle\bigwedge_{i \leq j}^{m}\bigwedge_{k = 1}^{\beta}ds_{ij}^{(k)}, &  \\
                      \displaystyle\bigwedge_{i=1}^{m} ds_{ii}\bigwedge_{i < j}^{m}\bigwedge_{k = 1}^{\beta}ds_{ij}^{(k)}, &
                       \hbox{if } s_{ii} \in \Re.
                    \end{array}
                  \right.
$$
\begin{rem}
Since generally the context establishes the conditions on the elements of
$\mathbf{S}$, that is, if $s_{ij} \in \Re$, $\in \mathfrak{C}$, $\in \mathfrak{H}$ or
$ \in \mathbf{O}$. It is considered that
$$
  (d\mathbf{S}) = \bigwedge_{i \leq j}^{m}\bigwedge_{k = 1}^{\beta}ds_{ij}^{(k)}
   \equiv \bigwedge_{i=1}^{m} ds_{ii}\bigwedge_{i < j}^{m}\bigwedge_{k =
1}^{\beta}ds_{ij}^{(k)}.
$$
Observe, too, that for the Lebesgue measure $(d\mathbf{S})$ defined thus, it is
required that $\mathbf{S} \in \mathfrak{P}_{m}^{\beta}$, that is, $\mathbf{S}$ must
be a non singular Hermitian matrix (Hermitian definite positive matrix). In the real
case, when $\mathbf{S}$ is a positive semidefinite matrix its corresponding measure
is studied in  \citet{u:94}, \citet{dg:97}, \citet{dggf:04a} and \citet{dggf:04b}
under different coordinate systems.
\end{rem}
 If $\mathbf{\Lambda} \in \mathfrak{D}_{m}^{\beta}$ then
$(d\mathbf{\Lambda})$ (the Legesgue measure in $\mathfrak{D}_{m}^{\beta}$) denotes
the exterior product of the $\beta m$ functionally independent variables
$$
  (d\mathbf{\Lambda}) = \bigwedge_{i = 1}^{n}\bigwedge_{k = 1}^{\beta}d\lambda_{i}^{(k)}.
$$
If $\mathbf{H}_{1} \in \mathcal{V}_{m,n}^{\beta}$ then
$$
  (\mathbf{H}^{*}_{1}d\mathbf{H}_{1}) = \bigwedge_{i=1}^{n} \bigwedge_{j =i+1}^{m}
  \mathbf{h}_{j}^{*}d\mathbf{h}_{i}.
$$
where $\mathbf{H} = (\mathbf{H}_{1}|\mathbf{H}_{2}) = (\mathbf{h}_{1}, \dots,
\mathbf{h}_{m}|\mathbf{h}_{m+1}, \dots, \mathbf{h}_{n}) \in \mathfrak{U}^{\beta}(m)$.
It can be proved that this differential form does not depend on the choice of the
matrix $\mathbf{H}_{2}$ and that it is invariant under the transformations
\begin{equation}\label{unif}
    \mathbf{H}_{1} \rightarrow \mathbf{QHP}_{1},  \quad \mathbf{Q }\in
    \mathfrak{U}^{\beta}(n) \hbox{ and } \mathbf{P} \in \mathfrak{U}^{\beta}(m).
\end{equation}
When $m = 1$; $\mathcal{V}^{\beta}_{1,n}$ defines the unit sphere in
$\mathfrak{F}^{n}$. This is, of course, an $(n-1)\beta$- dimensional surface in
$\mathfrak{F}^{n}$. When $m = n$ and denoting $\mathbf{H}_{1}$ by $\mathbf{H}$,
$(\mathbf{H}^{*}d\mathbf{H})$ is termed the \emph{Haar measure} on
$\mathfrak{U}^{\beta}(m)$ and defines an invariant differential form of a unique
measure $\nu$ on $\mathfrak{U}^{\beta}(m)$ given by
$$
  \nu (\mathfrak{M}) = \int_{\mathfrak{M}}(\mathbf{H}^{*}d\mathbf{H}).
$$
It is unique in the sense that any other invariant measure on
$\mathfrak{U}^{\beta}(m)$ is a finite multiple of $\nu$ and invariant because is
invariant under left and right translations, that is
$$
  \nu(\mathbf{Q}\mathfrak{M}) = \nu(\mathfrak{M}\mathbf{Q}) = \nu(\mathfrak{M}),
  \quad \forall \mathbf{Q} \in \mathfrak{U}^{\beta}(m).
$$
The surface area or volume of the Stiefel manifold $\mathcal{V}^{\beta}_{m,n}$ is
\begin{equation}\label{vol}
    \Vol(\mathcal{V}^{\beta}_{m,n}) = \int_{\mathbf{H}_{1} \in
  \mathcal{V}^{\beta}_{m,n}} (\mathbf{H}^{*}_{1}d\mathbf{H}_{1}) =
  \frac{2^{m}\pi^{mn\beta/2}}{\Gamma^{\beta}_{m}[n\beta/2]},
\end{equation}
where $\Gamma^{\beta}_{m}[a]$ denotes the multivariate Gamma function for the space
$\mathfrak{S}^{\beta}(m)$, and is defined by
\begin{eqnarray*}
  \Gamma_{m}^{\beta}[a] &=& \displaystyle\int_{\mathbf{A} \in \mathfrak{P}_{m}^{\beta}}
  \etr\{-\mathbf{A}\} |\mathbf{A}|^{a-(m-1)\beta/2 - 1}(d\mathbf{A}) \\
    &=& \pi^{m(m-1)\beta/4}\displaystyle\prod_{i=1}^{m} \Gamma[a-(i-1)\beta/2],
\end{eqnarray*}
where $\etr(\cdot) = \exp(\tr(\cdot))$, $|\cdot|$ denotes the determinant and $\re(a)
\geq (m-1)\beta/2$, see \citet{gr:87}.

\citet[p. 35]{d:02} first assigned the responsibility of the calculating some
Jacobianos to researchers of random matrix theory, although in the end this fact is
not especially important, as the Jacobian of the spectral factorisation (and of SVD
among others) for the symmetric matrix was first computed by \citet{j:54} and not by
Wigner (1958).

There exist an extensive bibliography on the computation of Jacobians of matrix
transformations and decompositions in both random matrix theory and statistics. We
now summarised diverse Jacobians in terms of the parameter $\beta$, some based on the
work of Dimitriu \citet{d:02}, while other results are proposed as extensions of real
or complex cases, see  \citet{do:51}, \citet{o:53}, \citet{ro:57}, \citet{j:54},
\citet{a:58}, \citet{j:64}, \citet{k:65}, \citet{m:82}, \citet{me:91}, \citet{o:02},
\citet{rva:05a} and \citet{rva:05b}. An excellent reference for real and complex
cases is \citet{m:97}, which includes many of the Jacobian that have been published.
Some Jacobians in the quaternionic case are obtained in \citet{lx:09}. We also
included a parameter count (or number of functionally independent variables, \#fiv),
that is, if $\mathbf{A}$ is factorised as $\mathbf{A} = \mathbf{BC}$, then the
parameter count is written as \#fiv in \textbf{A} = [\#fiv in \textbf{B}] + [\#fiv in
\textbf{C}], see \citet{d:02}.

\begin{lem}\label{lemlt}
Let $\mathbf{X}$ and $\mathbf{Y} \in {\mathcal L}_{m,n}^{\beta}$, and let $\mathbf{Y}
= \mathbf{AXB} + \mathbf{C}$, where $\mathbf{A} \in {\mathcal L}_{n,n}^{\beta}$,
$\mathbf{B} \in {\mathcal L}_{m,m}^{\beta}$ and $\mathbf{C} \in {\mathcal
L}_{m,n}^{\beta}$ are constant matrices. Then
\begin{equation}\label{lt}
    (d\mathbf{Y}) = |\mathbf{A}^{*}\mathbf{A}|^{\beta m/2} |\mathbf{B}^{*}\mathbf{B}|^{\beta
    n/2}(d\mathbf{X}).
\end{equation}
\end{lem}

\begin{lem}\label{lemhlt}
Let $\mathbf{X}$ and $\mathbf{Y} \in \mathfrak{P}_{m}^{\beta}$, and let $\mathbf{Y} =
\mathbf{AXA^{*}} + \mathbf{C}$, where $\mathbf{A}$ and $\mathbf{C} \in {\mathcal
L}_{m,m}^{\beta}$ are constant matrices. Then
\begin{equation}\label{hlt}
    (d\mathbf{Y}) = |\mathbf{A}^{*}\mathbf{A}|^{\beta(m-1)/2+1} (d\mathbf{X}).
\end{equation}
\end{lem}

\begin{lem}\label{lemdt}
Let $\mathbf{J} \in \mathfrak{T}_{U}^{\beta}(m)$ and write $\mathbf{J} = \mathbf{BG}$
where $\mathbf{B} = \diag(b_{1}, \dots, b_{m}) \in \mathfrak{D}_{m}^{\beta}$, and
$\mathbf{G} \in \mathfrak{T}_{U}^{\beta}(m)$, with $g_{ii}=1$ for all $i = 1, \dots,
m$. Then

$\bullet$ parameter count: $\beta m(m+1)/2 = [\beta m] + [\beta m(m-1)/2]$ and
\begin{equation}
     (d\mathbf{J}) = \prod_{i = 1}^{m} |b_{i}|^{\beta(m-i)} (d\mathbf{B})(d\mathbf{G}),
\end{equation}
where
$$
  (d\mathbf{G}) = \bigwedge_{i < j}^{m}\bigwedge_{k = 1}^{\beta}  dg_{ij}^{(k)}.
$$
\end{lem}

\begin{lem} [$LU$ decomposition, Doolittle's version]\label{lemdo}
Let $\mathbf{X} \in {\mathcal L}_{m,m}^{\beta}$, and write $\mathbf{X} =
\mathbf{\Delta}\mathbf{\Upsilon}$, where $\mathbf{\Delta} \in
\mathfrak{T}_{L}^{\beta}(m)$, with $\delta_{ii} = 1$, $i =1, \dots m$ and
$\mathbf{\Upsilon} \in \mathfrak{T}_{U}^{\beta}(m)$. Then

$\bullet$ parameter count: $\beta m^{2} = [\beta m(m-1)/2] + [\beta m(m+1)/2]$ and
\begin{equation}
    (d\mathbf{X}) = \prod_{i = 1}^{m} |\upsilon_{ii}|^{\beta(m-i)} (d\mathbf{\Upsilon})
    (d\mathbf{\Delta}),
\end{equation}
where $(d\mathbf{\Delta})= \bigwedge_{i > j}^{m} \bigwedge_{k = 1}^{\beta}
d\delta_{ij}^{(k)}$.
\end{lem}

\begin{lem}[$LU$ decomposition, Crout's version]\label{lemcr}
Let $\mathbf{X} \in {\mathcal L}_{m,m}^{\beta}$, such that $\mathbf{X} =
\mathbf{\Delta}\mathbf{\Upsilon}$, where $\mathbf{\Delta} \in
\mathfrak{T}_{L}^{\beta}(m)$ and $\mathbf{\Upsilon} \in \mathfrak{T}_{U}^{\beta}(m)$
such that $\upsilon_{ii} = 1$ for all $i = 1, \dots, m$. Then

$\bullet$ parameter count: $\beta m^{2} = [\beta m(m+1)/2] + [\beta m(m-1)/2]$ and
\begin{equation}
    (d\mathbf{X}) = \prod_{i = 1}^{m} |\delta_{ii}|^{\beta(m-i)} (d\mathbf{\Upsilon})
    (d\mathbf{\Delta}),
\end{equation}
where $(d\mathbf{\Upsilon}) = \bigwedge_{i < j}^{m}\bigwedge_{k = 1}^{\beta}
d\upsilon_{ij}^{\beta}$.
\end{lem}

As a consequence of Lemmas \ref{lemdo} and \ref{lemdt} there follows

\begin{lem}[ $LDM$ decomposition]\label{lemldm}
Let $\mathbf{X} \in {\mathcal L}_{m,m}^{\beta}$, such that $\mathbf{X} =
\mathbf{\Psi}\mathbf{\Pi} \mathbf{\Xi}$, where $\mathbf{\Psi} \in
\mathfrak{T}_{L}^{\beta}(m)$ with $\psi_{ii} = 1$ for all $i = 1, \dots, m$,
$\mathbf{\Pi}= \diag(\pi_{1}, \dots, \pi_{m}) \in \mathfrak{D}_{m}^{\beta}$,  and
$\mathbf{\Xi} \in \mathfrak{T}_{U}^{\beta}(m)$  with $\xi_{ii} = 1$ for all $i = 1,
\dots, m$. Then

$\bullet$ parameter count: $\beta m^{2} = [\beta m(m-1)/2] + [\beta m] + [\beta
m(m-1)/2]$ and
\begin{equation}\label{ldm}
    (d\mathbf{X}) = \prod_{i=1}^{m}|\pi_{ii}|^{2\beta(m-i)}(d\mathbf{\Psi})(d\mathbf{\Pi})(d\mathbf{\Xi}).
\end{equation}
\end{lem}

\begin{lem}[$QR$ decomposition]\label{lemqr}
Let $\mathbf{X} \in {\mathcal L}_{m,n}^{\beta}$, there then exists an $\mathbf{H}_{1}
\in \mathcal{V}_{m,n}^{\beta}$ and a  $\mathbf{T} \in \mathfrak{T}_{U}^{\beta}(m)$
with real $t_{ii} > 0$, $i = 1, 2, \ldots , q$ such that $\mathbf{X} =
\mathbf{H}_{1}\mathbf{T}$. Then

$\bullet$ parameter count: $\beta mn = [\beta mn - \beta m(m-1)/2 - m] + [\beta
m(m-1)/2 + m]$ and
\begin{equation}\label{qr}
   (d\mathbf{X}) = \prod_{i = 1}^{m} t_{ii}^{\beta(n-i+1)-1} (\mathbf{H}_{1}^{*}d\mathbf{H}_{1})(d\mathbf{T}).
\end{equation}
\end{lem}

Now, from Lemmas \ref{lemqr} and \ref{lemdo} it follows that

\begin{lem}[Modified $QR$ decomposition ($QDR$)]\label{lemmqr}
Let $\mathbf{X} \in {\mathcal L}_{m,n}^{\beta}$, then there exist $\mathbf{H}_{1} \in
\mathcal{V}_{m,n}^{\beta}$, a diagonal matrix $\mathbf{N} = \diag(n_{1}, \dots,
n_{m}) \in \mathfrak{D}_{m}^{1}$, with $n_{1} > \cdots > n_{m}> 0$ and
$\mathbf{\Omega} \in \mathfrak{T}_{U}^{\beta}(m)$ with $\omega_{ii} = 1$, $i = 1, 2,
\ldots , m$ such that $\mathbf{X} = \mathbf{H}_{1}\mathbf{N\Omega}$. Then

$\bullet$ parameter count: $\beta mn = [\beta mn - \beta m(m-1)/2 - m] + [m] + [\beta
m(m-1)/2]$ and
\begin{equation}\label{qdr}
    (d\mathbf{X}) = 2^{-m} \prod_{i = 1}^{m} n_{i}^{\beta(n+m-2i+1)-1}(\mathbf{H}_{1}^{*}
   d\mathbf{H}_{1})(d\mathbf{N})(d\mathbf{\Omega}).
\end{equation}
\end{lem}

\begin{lem}[Polar decomposition]\label{lempd}
Let $\mathbf{X} \in {\mathcal L}_{m,n}^{\beta}$ and write $\mathbf{X} = \mathbf{P}_{1}
\mathbf{R}$, with $\mathbf{P}_{1} \in \mathcal{V}_{m,n}^{\beta}$, and $\mathbf{R}
\in \mathfrak{P}_{m}^{\beta}$. Then

$\bullet$ parameter count: $\beta mn = [\beta mn - \beta m(m-1)/2 - m] + [\beta
m(m-1)/2 + m]$ and
\begin{equation}\label{p}
    (d\mathbf{X}) = |\mathbf{R}|^{\beta(n - m +1)-1}  \prod_{i<j}^{m}(d_{i} + d_{j})^{\beta}
    (d\mathbf{R})(\mathbf{P}_{1}^{*}d\mathbf{P}_{1}),
\end{equation}
where $\mathbf{R} = \mathbf{QDQ}^{*}$ is the spectral decomposition of $\mathbf{R}$,
$\mathbf{Q} \in \mathfrak{U}^{\beta}(m)$ and $\mathbf{D} = \diag(d_{1}, \ldots,
d_{m}) \in \mathfrak{D}_{m}^{1}$ with $d_{1}> \cdots > d_{m} > 0$.
\end{lem}

\begin{lem}[Singular value decomposition, $SVD$]\label{lemsvd}
Let $\mathbf{X} \in {\mathcal L}_{m,n}^{\beta}$, such that $\mathbf{X} =
\mathbf{V}_{1}\mathbf{DW}^{*}$ with $\mathbf{V}_{1} \in {\mathcal V}_{m,n}^{\beta}$,
$\mathbf{W} \in \mathfrak{U}^{\beta}(m)$ and $\mathbf{D} = \diag(d_{1}, \cdots,d_{m})
\in \mathfrak{D}_{m}^{1}$, $d_{1}> \cdots > d_{m} > 0$. Then

$\bullet$ parameter count: $\beta mn = [\beta mn - \beta m(m-1)/2 - m -(\beta -1)m] +
[m] + [\beta m(m-1)/2 + m]$ and
\begin{equation}\label{svd}
    (d\mathbf{X}) = 2^{-m}\pi^{\tau} \prod_{i = 1}^{m} d_{i}^{\beta(n - m + 1) -1}
    \prod_{i < j}^{m}(d_{i}^{2} - d_{j}^{2})^{\beta} (d\mathbf{D}) (\mathbf{V}_{1}^{*}d\mathbf{V}_{1})
    (\mathbf{W}^{*}d\mathbf{W}),
\end{equation}
where
$$
  \tau = \left\{
             \begin{array}{rl}
               0, & \beta = 1; \\
               -m, & \beta = 2; \\
               -2m, & \beta = 4; \\
               -4m, & \beta = 8.
             \end{array}
           \right.
$$
\end{lem}

\begin{lem} [Cholesky's decomposition]\label{lemch}
Let $\mathbf{S} \in \mathfrak{P}_{m}^{\beta}$ and write $\mathbf{S} =
\mathbf{T}^{*}\mathbf{T}$, where $\mathbf{T} \in \mathfrak{T}_{U}^{\beta}(m)$ with
$t_{ii} > 0$, $i = 1, 2, \ldots , m$. Then

$\bullet$ parameter count: $\beta m(m-1)/2 + m = \beta m(m-1)/2 + m$ and
\begin{equation}\label{ch}
    (d\mathbf{S}) = 2^{m} \prod_{i = 1}^{m} t_{ii}^{\beta(m - i) + 1} (d\mathbf{T}).
\end{equation}
\end{lem}

From Lemmas \ref{lemch} and \ref{lemdo} it follows that

\begin{lem} [ $L^{*}DL$ decomposition]\label{lemldl}
Let $\mathbf{S} \in \mathfrak{P}_{m}^{\beta}$ and write $S =
\mathbf{\Omega}^{*}\mathbf{O} \mathbf{\Omega}$, where $\mathbf{\Omega} \in
\mathfrak{T}_{U}^{\beta}(m)$ with $\omega_{ii} = 1$, $i = 1, 2, \ldots, m$ and a
diagonal matrix $\mathbf{O} = \diag(o_{1}, \dots, o_{m}) \in \mathfrak{D}_{m}^{1}$,
with $o_{1} > \cdots > o_{m}>0$. Then

$\bullet$ parameter count: $\beta m(m-1)/2 + m = [\beta m(m-1)/2] + [m]$ and
\begin{equation}\label{ldl}
    (d\mathbf{S}) = \prod_{i = 1}^{m} o_{i}^{\beta(m - i)}(d\mathbf{\Omega})(d\mathbf{O}).
\end{equation}
\end{lem}

\begin{lem}[Hermitian positive definite square root]\label{lemsr}
Let $\mathbf{S}$ and $\mathbf{R} \in \mathfrak{P}_{m}^{\beta}$ such that $\mathbf{S}
= \mathbf{R}^{2}$. Then

$\bullet$ parameter count: $\beta m(m-1)/2 + m = \beta m(m-1)/2 + m$ and
\begin{equation}\label{sr}
    (d\mathbf{S}) = 2^{m} |\mathbf{R}| \prod_{i \leq j}^{q}(d_{i} + d_{j})^{\beta}
    (d\mathbf{R}),
\end{equation}
where $\mathbf{R} = \mathbf{QDQ}^{*}$ is the spectral decomposition of $\mathbf{R}$,
$\mathbf{Q} \in \mathfrak{U}^{\beta}(m)$ and $\mathbf{D} = \diag(d_{1}, \ldots,
d_{m}) \in \mathfrak{D}_{m}^{1}$ with $d_{1}> \cdots > d_{m} > 0$.
\end{lem}

\begin{lem}[ Spectral decomposition]\label{lemsd}
Let $\mathbf{S} \in \mathfrak{P}_{m}^{\beta}$. Then the spectral decomposition can be
written as $\mathbf{S} = \mathbf{W}\mathbf{\Lambda W}^{*}$, where $\mathbf{W} \in
\mathfrak{U}^{\beta}(m)$ and $\mathbf{\Lambda} = \diag(\lambda_{1}, \dots,
\lambda_{m}) \in \mathfrak{D}_{m}^{1}$, with $\lambda_{1}> \cdots> \lambda_{m}>0$.
Then

$\bullet$ parameter count: $\beta m(m-1)/2 + m = [\beta m(m+1)/2 - m -(\beta -1)m] +
[m]$ and
\begin{equation}\label{sd}
    (d\mathbf{S}) = 2^{-m} \pi^{\tau} \prod_{i < j}^{m} (\lambda_{i} - \lambda_{j})^{\beta}
    (d\mathbf{\Lambda})(\mathbf{W}^{*}d\mathbf{W}),
\end{equation}
where $\tau$ is given in Lemma \ref{lemsvd}.
\end{lem}

Finally, by combining Lemmas \ref{lemsvd}, \ref{lempd}, \ref{lemqr} and \ref{lemmqr}
with Lemmas \ref{lemsd}, \ref{lemsr}, \ref{lemch}, \ref{lemldl}, respectively, the
following result is obtained.

\begin{lem}\label{lemW}
Let $\mathbf{X} \in {\mathcal L}_{m,n}^{\beta}$, and write $\mathbf{X} =
\mathbf{V}_{1}\mathbf{D}\mathbf{W}^{*}$ ($SVD$), $\mathbf{X} =
\mathbf{V}_{1}\mathbf{R}$ (Polar decomposition), $\mathbf{X} =
\mathbf{V}_{1}\mathbf{T}$ (QR decomposition) or $\mathbf{X} =
\mathbf{V}_{1}\mathbf{N\Omega}$ (modified QR decomposition) and let $\mathbf{S} =
\mathbf{X}^{*}\mathbf{X} \in \mathfrak{P}_{m}^{\beta}.$ Then
\begin{equation}\label{w}
    (d\mathbf{X}) = 2^{-m} |\mathbf{S}|^{\beta(n - m + 1)/2 - 1}
    (d\mathbf{S})(\mathbf{V}_{1}^{*}d\mathbf{V}_{1}).
\end{equation}
Note that the Lebesgue measure $(d\mathbf{S})$ in (\ref{w}), can be factorised in
terms of Cholesky, $L^{*}DL$, Hermitian positive definite square root or spectral
decompositions, thus obtaining this way alternative explicit expressions of
$(d\mathbf{S})$, given by (\ref{ch}), (\ref{ldl}), (\ref{sr}) and (\ref{sd}),
respectively. For the corresponding coordinates and bases, see \citet[Chapter
6]{ea:07}.
\end{lem}

\section{Elliptical ensemble}\label{sec3}

In this section we employ the nomenclature recently proposed in random matrix theory
for some of the most commonly-studied random matrices. These terms are
\emph{Hermite}, \emph{Laguerre}, \emph{Jacobi} and \emph{Fourier} instead of
Gaussian, Wishart, \textsc{manova} (or beta type I) and circular, see \citet{de:02},
\citet{er:05}, \citet{es:08}. These alternative names are proposed in view of the
fact that in random matrix theory we are interested in matrices with joint eigenvalue
density proportional to $\prod_{i=1}^{m} w(\lambda_{i})|\Delta(\Lambda)|^{\beta}$
where $|\Delta(\Lambda)| = \prod_{i < j}(\lambda_{i}-\lambda_{j})$ is the absolute
value of the Vandermonde determinant and $w(\lambda)$ is a weight function of a
system of orthogonal polynomials. In particular, when the joint eigenvalue density is
calculated under Gaussian, Wishart and \textsc{manova} random matrices, these weight
functions correspond to the Hermite, Laguerre and Jacobi orthogonal polynomials,
respectively.

Under generalised multivariate analysis, four classes of matrix variate elliptical
distributions have been defined and studied, see \citet{fz:90}. A random matrix
$\mathbf{X} \in \mathcal{L}_{m,n}^{1}$ is said to have a matrix variate
left-elliptical distribution, the largest of the four classes of matrix variate
elliptical distributions, if its density function is given by
$$
  \frac{c(m,n)}{|\mathbf{\Sigma}|^{n/2}|\mathbf{\Theta}|^{m/2}}h\left (\mathbf{\Sigma}^{-1/2}
  (\mathbf{X} - \boldsymbol{\mu})^{T} \mathbf{\Theta}^{-1}(\mathbf{X} - \boldsymbol{\mu})
  \mathbf{\Sigma}^{-1/2}\right ),
$$
where $h$ is a real function, $c(m,n)$ denotes the normalisation constant which might
be a function of another parameter implicit in the function $h$, $\mathbf{\Sigma} \in
\mathfrak{P}_{m}^{\beta}$, $\mathbf{\Theta} \in \mathfrak{P}_{n}^{\beta}$ and
$\boldsymbol{\mu} \in \Re^{n \times m}$. This fact is denoted as
$$
  \mathbf{X} \sim \mathcal{ELS}_{m \times n}(\boldsymbol{\mu},\mathbf{\Sigma}, \mathbf{\Theta}, h).
$$
One might wish to consider situations where $\mathbf{X}$ has a density function of
the following form, see \citet{fz:90} and \citet{fl:99},
\begin{equation}\label{eliptica1}
    \frac{c(m,n)}{|\mathbf{\Sigma}|^{n/2}|\mathbf{\Theta}|^{m/2}}h\left (\mathbf{\Sigma}^{-1}
    (\mathbf{X} - \boldsymbol{\mu})^{T} \mathbf{\Theta}^{-1}(\mathbf{X} - \boldsymbol{\mu})\right).
\end{equation}
This condition is equivalent to considering the function $h$ as a symmetric function,
i.e. $g| g(\mathbf{AB}) = g(\mathbf{BA})$ for any symmetric matrices $\mathbf{A}$ and
$\mathbf{B}$. This condition is equivalent to that in which $h(\mathbf{\mathbf{A}})$
depends on $\mathbf{A}$ only through its eigenvalues, in which case the function
$h(\mathbf{A})$ can be expressed as $h(\lambda(\mathbf{A}))$, where
$\lambda(\mathbf{A})= \diag(\lambda_{1}, \dots, \lambda_{m})$ and $\lambda_{1},
\dots, \lambda_{m}$ are the eigenvalues of $\mathbf{A}$. Two subclasses of matrix
variate elliptical distributions are of particular interest: \emph{matrix variate
vector-spherical and spherical elliptical distributions}. For these distributions,
$\lambda (\mathbf{A}) \equiv \tr (\mathbf{A})$ and $\lambda(\mathbf{A})$ represent
any function of eigenvalues of $\mathbf{A}$, and are denoted as $\mathbf{X} \sim
\mathcal{EVS}_{m \times n}(\boldsymbol{\mu},\mathbf{\Sigma}, \mathbf{\Theta}, h)$ and
$\mathbf{X} \sim \mathcal{ESS}_{m \times n}(\boldsymbol{\mu},\mathbf{\Sigma},
\mathbf{\Theta}, h)$, respectively. Note that matrix variate vector-spherical
elliptical distributions are a subclass of matrix variate spherical elliptical
distributions. Many well-known distributions are examples of these subclasses; one
such is the matrix variate \emph{Hermite distribution}. Other variants include matrix
variate vector-spherical elliptical distributions, e.g. Pearson type II, Pearson type
VII, Kotz type , Bessel and Logistic, among many others, see \citet{gv:93}. Further
cases are those of matrix variate spherical elliptical distributions, e.g. Pearson
Type II, Pearson type VII and Kotz type, among many others, see \citet{fl:99}.

When in matrix variate vector-spherical and spherical elliptical distributions it is
assumed that $\mathbf{\Sigma} = \mathbf{I}_{m}$, $\mathbf{\Theta} = \mathbf{I}_{n}$
and $\boldsymbol{\mu}= \mathbf{0}_{n\times m}$ we obtain the matrix variate
vector-spherical (denoted as $\mathcal{VS}_{m \times n}(\mathbf{0},\mathbf{I}_{m},
\mathbf{I}_{n}, h)\equiv \mathcal{VS}_{m \times n}(h)$) and spherical distributions
(denoted as $\mathcal{SS}_{m \times n}(\mathbf{0},\mathbf{I}_{m}, \mathbf{I}_{n},
h)\equiv \mathcal{SS}_{m \times n}(h)$), which are known to be \emph{invariant
distributions under orthogonal transformations}, because $\mathbf{X}$ and
$\mathbf{QXP}$ have the same distribution when $\mathbf{P} \in \mathfrak{U}^{1}(m)$
and $\mathbf{Q} \in \mathfrak{U}^{1}(n)$. In particular, note that if $x$ is a random
variable with vector-spherical or spherical distribution it is denoted as $x \sim
\mathcal{VS}(0,1,h)\equiv \mathcal{VS}(h)$ or $x \sim \mathcal{SS}(0,1,h)\equiv
\mathcal{SS}(h)$, respectively.

We now define \emph{spherical} and \emph{generalised Laguerre ensembles} based on
vector-spherical and spherical distributions. As can be seen, these ensembles contain
as particular cases the classical Hermite, Laguerre, Jacobi and Fourier ensembles, as
well as many others that have been studied in the literature of random matrix theory,
see \citet{f:05}.

\subsection{Vector-spherical and spherical random matrices}\label{subsec31}

$VS^{\beta}(n,m,h)$  and $SS^{\beta}(n,m,h)$ are $n \times m$ matrices of
non-correlated and identically distributed (n-c.i.d.) with entries
$\mathcal{VS}^{\beta}(0,1,h)$  and $\mathcal{SS}^{\beta}(0,1,h)$, respectively.

\begin{enumerate}
\item If $\mathbf{A} \in \mathcal{L}_{m,n}^{\beta}$ is an $n \times m$ vector-spherical
random matrix $VS^{\beta}(n,m,h)$ then its joint element density is given by
\begin{equation}\label{vsrm}
    c^{\beta}(m,n)h(\tr(\mathbf{A}^{*}\mathbf{A}))
\end{equation}
where
$$
  c^{\beta}(m,n) = \frac{\Gamma[\beta mn/2]}{2 \pi^{\beta mn/2}}\left\{\int_{\mathfrak{P}_{1}^{\beta}}
    u^{\beta mn -1}h(u^{2}) du \right\}^{-1},
$$
\item If $\mathbf{A} \in \mathcal{L}_{m,n}^{\beta}$ is an $n \times m$ spherical random
matrix $SS^{\beta}(n,m,h)$ then its joint element density is given by
\begin{equation}\label{ssrm}
    d^{\beta}(m,n)h(\lambda(\mathbf{A}^{*}\mathbf{A}))
\end{equation}
where $d^{\beta}(m,n)$ is given by
\begin{equation}\label{dc}
    \frac{\Gamma_{m}^{\beta}[\beta n/2]\Gamma_{m}^{\beta}[\beta m/2]}{\pi^{m(n+1)\beta/2+\tau}}
    \left\{\build{\int \cdots \int}{}{\l_{1} > \cdots > l_{m} > 0} |\mathbf{L}|^{\beta(n-m+1)/2-1}
    |\Delta(\mathbf{L})|^{\beta}h(\mathbf{L}) (d\mathbf{L})\right\}^{-1},
\end{equation}
with $\mathbf{L} = \diag(l_{1}, \dots, l_{m})$ and $l_{i}$ is the $i$th eigenvalue of
$(\mathbf{A}^{*}\mathbf{A})$ and $\tau$ is given in Lemma \ref{lemsvd}, for the real
case see \citet{fkn:90}, \citet{fz:90} and \citet{gv:93}.
\end{enumerate}

Note that $\mathbf{A}$ and $\mathbf{QAP}$ have the same distribution when $\mathbf{P}
\in \mathfrak{U}^{\beta}(m)$ and $\mathbf{Q} \in \mathfrak{U}^{\beta}(n)$. That is,
the distributions of $VS^{\beta}(n,m,h)$ and $SS^{\beta}(n,m,h)$ \emph{are invariant}
under orthogonal ($\beta = 1$), unitary ($\beta = 2$), symplectic ($\beta = 4$) and
exceptional type ($\beta = 8$) transformations.

\begin{rem}
Note that the unique case in which the elements $a_{ij}$ of the random matrix
$\mathbf{A}$ are i.i.d. is when
$\mathcal{VS}^{\beta}(0,1,h)=\mathcal{SS}^{\beta}(0,1,h) = \mathcal{N}^{\beta}(0,1)$,
this is, when $a_{ij}$ has a standard normal (Gaussian or of Hermite) distribution,
see \citet[Theorem 2.7.1, p. 72]{fz:90} and \citet[Theorems 6.2.1 and 6.2.2
p.193]{gv:93}.
\end{rem}

Some particular vector-spherical random matrices and their corresponding density
functions are summarised in the following result, where $u \equiv \tr
\mathbf{A}^{*}\mathbf{A})$.

\begin{cor}\label{corvsrm}
Let $\mathbf{A} \in \mathcal{L}_{m,n}^{\beta}$ with density function (\ref{vsrm})
then
\begin{enumerate}
  \item \textbf{(random matrix of Hermite)} If $h(u) = \exp (-\beta u/2)$, it is said that
  $\mathbf{A}$ is a random matrix of Hermite with density function
  $$
    \frac{1}{\left(2 \pi \beta^{-1}\right)^{\beta mn/2}} \etr\{-\beta \mathbf{A}^{*}
    \mathbf{A}/2\}.
  $$
  \item $(\mathbf{T}$ \textbf{type I random matrix)} If $h(u) = (1+u/\nu)^{-\beta(mn + \nu/2)}$,
  it is said that $\mathbf{A}$ is a $T$ type I random matrix with density function
  $$
    \frac{\Gamma^{\beta}[\beta(mn +\nu)/2]}{\pi^{\beta mn/2} \Gamma^{\beta}[\beta \nu/2]}
    (1 + \tr \mathbf{A}^{*}\mathbf{A})^{-\beta(mn + \nu)/2},
  $$
  where $\nu >0$. This distribution is also termed the Pearson type VII matrix variate
  distribution, see \citet{gv:93} and \citet{dggj:09}.
  \item \textbf{(Gegenbauer type I random matrix)} If $h(u) = (1-u)^{-\beta q}$,
  it is said that $\mathbf{A}$ is a Gegenbauer type I random matrix with density function
  $$
    \frac{\Gamma^{\beta}[\beta mn/2 + \beta q +1]}{\pi^{\beta mn/2} \Gamma^{\beta}[\beta q +1]}
    (1 - \tr \mathbf{A}^{*}\mathbf{A})^{-\beta q},
  $$
  where $q > -1$ and $\tr \mathbf{A}^{*}\mathbf{A} \leq 1$.
  This distribution is known in statistical bibliography as the inverted $T$ or Pearson type II matrix
  variate distribution, see \citet{gv:93}, \citet{p:82} and \citet{dggj:09}.
\end{enumerate}
\end{cor}

For the spherical random matrices we have.

\begin{cor}\label{corssrm}
Let $\mathbf{A} \in \mathcal{L}_{m,n}^{\beta}$ with density function (\ref{vsrm})
then
\begin{enumerate}
  \item $(\mathbf{T}$ \textbf{type II random matrix)} If $h(\lambda(\mathbf{A}^{*}\mathbf{A})) =
  |\mathbf{I}_{m} + \mathbf{A}^{*}\mathbf{A}|^{-\beta(n + \nu)/2}$,
  it is said that $\mathbf{A}$ is a $T$ type II random matrix with density function
  $$
    \frac{\Gamma_{m}^{\beta}[\beta(n +\nu)/2]}{(\pi)^{\beta mn/2} \Gamma_{m}^{\beta}[\beta \nu/2]}
    |\mathbf{I}_{m} + \mathbf{A}^{*}\mathbf{A}|^{-\beta(n + \nu)/2},
  $$
  where $\nu > m$. This distribution is also termed the $T$ or Pearson type VII matricvariate
  distribution, see \citet{di:67}, \citet{p:82}, \citet{fl:99} and \citet{dggj:09}.
  \item \textbf{(Gegenbauer type II random matrix)} If
  $$
    h(\lambda(\mathbf{A}^{*}\mathbf{A})) = |I - \mathbf{A}^{*}\mathbf{A}|^{\beta(\nu - m +1)/2 -1},
  $$
  it is said that $\mathbf{A}$ is a Gegenbauer type II random matrix with density function
  $$
    \frac{\Gamma_{m}^{\beta}[\beta (n + \nu)/2]}{\pi^{\beta mn/2} \Gamma_{m}^{\beta}[\beta \nu/2]}
    |I - \mathbf{A}^{*}\mathbf{A}|^{\beta(\nu - m +1)/2 -1},
  $$
  where $\nu > (m-1)/2$ and $\mathbf{A}^{*}\mathbf{A} \in \mathfrak{P}_{m}^{\beta}$.
  This distribution is known in statistical bibliography as the inverted $T$ or Pearson type II
  matricvariate distribution, see \citet{p:82}, \citet{fl:99} and \citet{dggj:09}.
\end{enumerate}
\end{cor}

\subsection{Construction of the vector-spherical and spherical random matrix ensembles}\label{subsec32}

The vector-spherical and generalised vector-spherical-Laguerre ensembles are
constructed from $VS^{\beta}(n,m,h)$ as follows.

\textbf{Vector-spherical orthogonal ensemble (VSOE):} symmetric $m \times m$ matrix
obtained as $(\mathbf{A} + \mathbf{A}^{T})/2$, where $\mathbf{A}$ is $VS^{1}(m,m,h)$.
The diagonal entries are n-c.i.d. with distribution $\mathcal{VS}^{1}(0,1,h)$, and
the off-diagonal entries are n-c.i.d. (subject to the symmetry) with distribution
$\mathcal{VS}^{1}(0,1/2,h)$.

\textbf{Vector-spherical unitary ensemble (VSUE):} Hermitian $m \times m$ matrix
obtained as $(\mathbf{A} + \mathbf{A}^{*})/2$, where $\mathbf{A}$ is $VS^{2}(m,m,h)$
and $*$ denotes the Hermitian transpose of a complex matrix. The diagonal entries are
n-c.i.d. with distribution $\mathcal{VS}^{1}(0,1,h)$, and the off-diagonal entries
are n-c.i.d. (subject to the symmetry) with distribution $\mathcal{VS}^{2}(0,1/2,h)$.

\textbf{Vector-spherical symplectic ensemble (VSSE):} quaternionic Hermitian
(self-dual) $m \times m$ matrix obtained as $(\mathbf{A} + \mathbf{A}^{*})/2$, where
$\mathbf{A}$ is $VS^{4}(m,m,h)$ and $*$ denotes the quaternionic Hermitian (or dual)
transpose of a quaternion matrix. The diagonal entries are n-c.i.d. with distribution
$\mathcal{VS}^{1}(0,1,h)$, and the off-diagonal entries are n-c.i.d. (subject to the
symmetry) with distribution $\mathcal{VS}^{4}(0,1/2,h)$.

\textbf{Vector-spherical exceptional type ensemble (VSETE):} octonionic Hermitian $m
\times m$ matrix obtained as $(\mathbf{A} + \mathbf{A}^{*})/2$, where $\mathbf{A}$ is
$VS^{8}(m,m,h)$ and $*$ denotes the octonionic Hermitian transpose of a octonionic
matrix. The diagonal entries are n-c.i.d. with distribution
$\mathcal{VS}^{1}(0,1,h)$, and the off-diagonal entries are n-c.i.d. (subject to the
symmetry) with distribution $\mathcal{VS}^{8}(0,1/2,h)$.

Analogously, by replacing $\mathcal{VS}$ and $LS$ by $\mathcal{SS}$ and $SS$,
respectively, in the four previous definitions we obtain the \textbf{Spherical
orthogonal ensemble (SSOE)}, \textbf{Spherical unitary ensemble (SSUE)},
\textbf{Spherical symplectic ensemble (SSSE)}, and \textbf{Spherical exceptional type
ensemble (SSETE)}.

Similarly, generalised vector-spherical-Laguerre and spherical-Laguerre ensembles can
be defined as follows.

\textbf{Generalised vector-spherical-Laguerre} $\hbox{\boldmath($VSL^{\beta}(n,m,h)$,
$n \geq m$):}$ Symmetric/ Hermitian/ quaternionic Hermitian/ octonionic Hermitian $m
\times m$ matrix which can be obtained as $\mathbf{A}^{*}\mathbf{A}$, where
$\mathbf{A}$ is $VS^{\beta}(n,m,h)$.

Again, by replacing $VS^{\beta}$ by $SS^{\beta}$ in the previous definition we obtain
the \textbf{Generalised spherical-Laguerre} $\hbox{\boldmath($SL^{\beta}(n,m,h)$, $n
\geq m$)}$.

As shown later, the Jacobi ensemble is a particular case of $VSL^{\beta}(n,m,h)$ or
$SL^{\beta}(n,m,h)$ and the Fourier ensemble is a particular case of
$SS^{\beta}(m,m,h)$.

\subsection{Computing the joint element densities}\label{subsec33}

Let $\mathbf{A}$ be an $m \times m$ matrix from the vector-spherical ensemble, that
is
\begin{equation}\label{vs1}
    a_{ij} \sim \left\{
                \begin{array}{ll}
                  \mathcal{VS}^{\beta}(0,1,h), & i = j, \\
                  \mathcal{VS}^{\beta}(0,1/2,h), & i > j.
                \end{array}
              \right.
\end{equation}
Then, it is straightforward to see that its joint element density is
\begin{equation}\label{vsed}
    2^{m(m-1)\beta/4}c^{\beta}(m,m)h\left(\tr \mathbf{A}^{2}\right),
\end{equation}
where
$$
  c^{\beta}(m,m) = \frac{\Gamma[(m + m(m - 1)\beta/2)/2]}{2 \pi^{(m + m(m - 1)\beta/2)/2}}
  \left\{\int_{\mathfrak{P}_{1}^{\beta}} u^{m + m(m - 1)\beta/2 - 1}h(u^{2}) du \right\}^{-1}.
$$
\begin{rem}
In the general case, from (\ref{eliptica1}) we have $|\mathbf{\Sigma}|^{\beta
n/2}|\mathbf{\Theta}|^{\beta m/2} = |\mathbf{\Theta} \otimes
\mathbf{\Sigma}|^{\beta/2}$. Under (\ref{vs1}), $\mathbf{\Theta} \otimes
\mathbf{\Sigma}$ is a diagonal matrix with $m$-times ones and $m(m-1)/2$-times $1/2$,
considering that $\mathbf{A}$ is symmetric, from which $|\mathbf{\Theta} \otimes
\mathbf{\Sigma}|^{\beta/2} = (1/2)^{m(m-1)\beta/4}$.
\end{rem}
Similarly, let $\mathbf{A}$ be an $m \times m$ matrix from the spherical ensemble,
that is
\begin{equation}\label{ss1}
    a_{ij} \sim \left\{
                \begin{array}{ll}
                  \mathcal{SS}^{\beta}(0,1,h), & i = j, \\
                  \mathcal{SS}^{\beta}(0,1/2,h), & i > j.
                \end{array}
              \right.
\end{equation}
And then it follows that its joint element density is
\begin{equation}\label{ssed}
    2^{m(m-1)\beta/4}d^{\beta}(m,m)h\left(\lambda\left(\mathbf{A}^{2}\right)\right),
\end{equation}
where
\begin{equation}\label{dc2}
    d^{\beta}(m,m) =\frac{\Gamma_{m}^{\beta}[\beta m/2]}{2^{m(m-1)\beta/4}\pi^{\beta m^{2}/2+\tau}}
    \left\{\build{\int \cdots \int}{}{\l_{1} > \cdots > l_{m} > 0}
    |\Delta(\mathbf{L})|^{\beta}h(\mathbf{L^{2}}) (d\mathbf{L}) \right\}^{-1},
\end{equation}
with $\mathbf{L} = \diag(l_{1}, \dots, l_{m})$, where $l_{i}$ is the $i$th eigenvalue
of $\mathbf{A}$ and $\tau$ is as given in Lemma \ref{lemsvd}.

Consider the generalised vector-spherical-Laguerre ensemble  $VL(n,m,h) = \mathbf{S}
= \mathbf{A}^{*}\mathbf{A}$, where $\mathbf{A} = VS(n,m,h)$. Its joint element
density can be computed in diverse ways. For example, following the approach of
\citet{h:55}, let $\mathbf{A} = \mathbf{V}_{1}\mathbf{R}$ be the polar decomposition
of $\mathbf{A}$, then $\mathbf{S} = \mathbf{A}^{*}\mathbf{A} = \mathbf{R}^{2}$, from
Lemma \ref{w} and (\ref{vsrm}), the joint density of $\mathbf{S}$ and
$\mathbf{V}_{1}$ is
$$
  c^{\beta}(m,n)h(\tr(\mathbf{S}))2^{-m} |\mathbf{S}|^{\beta(n - m + 1)/2 - 1}
    (d\mathbf{S})(\mathbf{V}_{1}^{*}d\mathbf{V}_{1}).
$$
On integrating with respect to $\mathbf{V}_{1}$ using (\ref{vol}) the marginal
density of $\mathbf{S}$ is found to be
\begin{equation}\label{vsld}
    \frac{c^{\beta}(m,n) \pi^{\beta mn/2}}{\Gamma_{m}^{\beta}[\beta n/2]} |\mathbf{S}|^{\beta(n - m + 1)/2 - 1}
    h(\tr(\mathbf{S})).
\end{equation}
Note that the same result can be obtained from Lemma \ref{w} and taking into account
the SV, QR or MQR decomposition instated of the polar decomposition in the previous
procedure.

As in the generalised vector-spherical-Laguerre ensemble case, it is obtained that
the joint element density of the generalised spherical-Laguerre ensemble is
\begin{equation}\label{sld}
    \frac{d^{\beta}(m,n) \pi^{\beta mn/2}}{\Gamma_{m}^{\beta}[\beta n/2]} |\mathbf{S}|^{\beta(n - m + 1)/2 - 1}
    h(\lambda(\mathbf{S})).
\end{equation}

As was done for Corollary \ref{corvsrm}, some particular vector-spherical random
matrix ensembles and their corresponding density functions are summarised in the
following result.

\begin{cor}\label{corvspd}
Let $\mathbf{A} \in \mathfrak{P}_{m}^{\beta}$ with density function (\ref{vsed}) then
\begin{enumerate}
  \item \textbf{(Classical Hermite ensemble)} If $h(u) = \exp (-\beta u/2)$, it is said that
  $\mathbf{A}$ is a Hermite ensemble with density function
  $$
    \frac{1}{2^{m/2}\left(\pi \beta^{-1}\right)^{m/2 + m(m-1)\beta/4}} \etr\{-\beta \mathbf{A}^{2}
    /2\}.
  $$
  \item $(\mathbf{T}$ \textbf{type I ensemble)} If $h(u) = (1+u/\nu)^{-m/2 - \beta(m(m-1)/2 + \nu)/2}$,
  it is said that $\mathbf{A}$ is a $T$ type I ensemble with density function
  $$
    \frac{2^{m(m-1)\beta/4}\Gamma^{m/2\beta}[\beta(m(m-1)/2 +\nu)/2]}{\pi^{m/2 + m(m-1)\beta/4} \Gamma^{\beta}[\beta \nu/2]}
    (1 + \tr \mathbf{A}^{2})^{-m/2 - \beta(m(m-1)/2 + \nu)/2},
  $$
  where $\nu >0$. This ensemble could be termed the Pearson type VII matrix variate
  ensemble.
  \item \textbf{(Gegenbauer type I ensemble)} If $h(u) = (1-u)^{-\beta q}$,
  it is said that $\mathbf{A}$ is a Gegenbauer type I ensemble with density function
  $$
    \frac{2^{m(m-1)\beta/4}\Gamma^{\beta}[m/2 +m(m-1)\beta/4 + \beta q +1]}{\pi^{m/2 + m(m-1)\beta /4}
    \Gamma^{\beta}[\beta q +1]} (1 - \tr \mathbf{A}^{2})^{-\beta q},
  $$
  where $q > -1$ and $\tr \mathbf{A}^{2} \leq 1$.
  In this case, too, the ensemble might be termed the inverted T or Pearson type II matrix variate ensemble.
\end{enumerate}
\end{cor}

For the spherical random matrix ensembles we have the following.

\begin{cor}\label{corsspd}
Let $\mathbf{A} \in \mathfrak{P}_{m}^{\beta}$ with density function (\ref{ssed}) then
\begin{enumerate}
  \item $(\mathbf{T}$ \textbf{type II ensemble)} If $h(\lambda(\mathbf{A}^{2})) =
  |\mathbf{I}_{m} + \mathbf{A}^{2}|^{-\beta(n + \nu)/2}$,
  it is said that $\mathbf{A}$ is a $T$ type II ensemble with density function
  $$
    2^{m(m-1)\beta/4} c^{\beta}(m,m)
    |\mathbf{I}_{m} + \mathbf{A}^{2}|^{-\beta(n + \nu)/2},
  $$
  where $\nu > m$. This ensemble could possibly be termed the Pearson type VII matricvariate
  ensemble.
  \item \textbf{(Gegenbauer type II ensemble)} If $h(\lambda(\mathbf{A}^{2})) = |\mathbf{I}_{m} -
  \mathbf{A}^{2}|^{\beta(\nu - m +1)/2 -1}$,
  it is said that $\mathbf{A}$ is a Gegenbauer type II ensemble with density function
  $$
     2^{m(m-1)\beta/4} c^{\beta}(m,m)
    |\mathbf{I}_{m} - \mathbf{A}^{2}|^{\beta(\nu - m +1)/2 -1},
  $$
  where $\nu \geq (m-1)\beta$ and $\mathbf{A}^{*}\mathbf{A} \in \mathfrak{P}_{m}^{\beta}$.
  This ensemble could be termed an inverted $T$ or Pearson type II matricvariate
  ensemble.
\end{enumerate}
\end{cor}

Analogously, some particular cases for vector-spherical and spherical-Laguerre random
matrix ensembles are obtained in the following two results.

\begin{cor}
Let $\mathbf{S} \in \mathfrak{P}_{m}^{\beta}$ with density function (\ref{vsld}) then
\begin{enumerate}
  \item \textbf{(Classical Laguerre ensemble)} If $h(u) = \exp (-\beta u/2)$, it is said that
  $\mathbf{S}$ is a Laguerre ensemble with density function
  $$
    \frac{1}{\left(2\beta^{-1}\right)^{\beta mn/2} \Gamma_{m}^{\beta}[\beta mn/2]}
    |\mathbf{S}|^{\beta(n - m +1)/2 -1}\etr\{-\beta \mathbf{S}/2\}
  $$
  where $n \geq (m-1)\beta$.
  \item $(\mathbf{T}$-\textbf{Laguerre type I ensemble)} If $h(u) = (1+u/\nu)^{-\beta(n + \nu)/2}$,
  it is said that $\mathbf{S}$ is a $T$-Laguerre type I ensemble with density function
  $$
    \frac{\Gamma^{\beta}[\beta(mn+\nu)/2]}{\Gamma^{\beta}[\beta \nu/2]
    \Gamma_{m}^{\beta}[\beta n/2]} |\mathbf{S}|^{\beta(n - m +1)/2 -1}
    (1 + \tr \mathbf{S})^{-\beta(n + \nu)/2},
  $$
  where $\nu >0$ and $n \geq (m-1)\beta$.
  \item \textbf{(Gegenbauer-Laguerre type I ensemble)} If $h(u) = (1-u)^{-\beta q}$,
  it is said that $\mathbf{S}$ is a Gegenbauer-Laguerre type I ensemble with density function
  $$
    \frac{\Gamma^{\beta}[\beta mn/2 + \beta q +1]}{ \Gamma^{\beta}[\beta q +1]
    \Gamma_{m}^{\beta}[\beta n/2]} |\mathbf{S}|^{\beta(n - m +1)/2 -1}
    (1 - \tr \mathbf{S})^{\beta q},
  $$
  where $q > -1$, $n > (m-1)\beta/2$ and $\tr \mathbf{S} \leq 1$.
\end{enumerate}
\end{cor}

For the spherical-Laguerre random matrix ensembles we have the following.

\begin{cor}
Let $\mathbf{S} \in \mathfrak{P}_{m}^{\beta}$ with density function (\ref{sld}) then
\begin{enumerate}
  \item $(\mathbf{T}$\textbf{-Laguerre type II ensemble)} If $h(\lambda(\mathbf{S})) =
  |\mathbf{I}_{m} + \mathbf{S}|^{-\beta(n + \nu)/2}$,
  it is said that $\mathbf{S}$ is a $T$-Laguerre type II ensemble with density function
  $$
    \frac{1}{\mathcal{B}_{m}^{\beta}[\beta m/2,\beta n/2]} |\mathbf{S}|^{\beta(n - m +1)/2 -1}
    |\mathbf{I}_{m} + \mathbf{S}|^{-\beta(n + \nu)/2},
  $$
  where $\nu \geq (m-1)\beta$ and $n \geq (m-1)\beta$. This distribution is
  also known as the Studentised Wishart distribution, see \citet{o:64}.
  \item \textbf{(Gegenbauer-Laguerre type II ensemble)} If $h(\lambda(\mathbf{S}))
  = |\mathbf{I}_{m} - \mathbf{S}|^{\beta(\nu - m +1)/2 -1}$,
  it is said that $\mathbf{S}$ is a Gegenbauer-Laguerre type II ensemble with density function
  $$
     \frac{1}{\mathcal{B}_{m}^{\beta}[\beta n/2,\beta \nu/2]} |\mathbf{S}|^{\beta(n - m +1)/2 -1}
    |\mathbf{I}_{m} - \mathbf{S}|^{\beta(\nu - m +1)/2 -1},
  $$
  where $\nu  \geq (m-1)\beta$ and $n \geq (m-1)\beta$.
\end{enumerate}
Finally, note that $\mathcal{B}_{m}^{\beta}[a,b]$, defined as
$$
  \mathcal{B}_{m}^{\beta}[a,b] =
\frac{\Gamma_{m}^{\beta}[a]\Gamma_{m}^{\beta}[a]}{\Gamma_{m}^{\beta}[a+b]}
$$
is the multivariate beta function, where $\re(a) > (m-1)\beta/2$ and $\re(b) >
(m-1)\beta/2$, see \citet{h:55}.
\end{cor}

\begin{rem}\label{rem33}
\begin{enumerate}
  \item Note that, from (\ref{unif}), $\mathbf{H}_{1} \in \mathcal{V}_{m,n}^{\beta}$ is a
    spherical random matrix. Moreover, (see \citet[Lemma 3.1.3(iii), p. 94]{fz:90}), the differential
    form of its density function is
    $$
      \frac{1}{\Vol\left(\mathcal{V}_{m,n}^{\beta}\right)}
      (\mathbf{H}_{1}^{*}d\mathbf{H}_{1}) = \frac{\Gamma_{m}^{\beta}[\beta n/2]}{2^{m}\pi^{\beta mn/2}}
      (\mathbf{H}_{1}^{*}d\mathbf{H}_{1}).
    $$
    In the context of random matrix theory, $\mathbf{H}_{1}$ is a Fourier random
    matrix. Now, let $\mathfrak{U}^{\beta}_{S}(m)$ be the group of orthogonal Hermitian
    matrices $\mathbf{H}$. Then $\mathbf{H}$ is a Fourier ensemble with a differential form of its density
    function
    \begin{equation}\label{vol2}
        \frac{1}{\Vol\left(\mathfrak{P}_{S}^{\beta}(m)\right)}
      (\mathbf{H}d\mathbf{H}) = \frac{\Gamma_{m}^{\beta}[\beta m/2] \Gamma[\beta m/2 +1]}{2^{m}
      \pi^{\beta m^{2}/2 + \tau} (\Gamma[\beta/2 +1])^{m}}(\mathbf{H}d\mathbf{H}).
    \end{equation}
    The value of $\Vol\left(\mathfrak{P}_{S}^{\beta}(m)\right)$ is found in next
    subsection.
  \item Some $g(\mathbf{A})$ functions, where $\mathbf{A}$ is a spherical random
    matrix, have invariant distributions under the corresponding class of spherical
    distribution, under certain conditions; in other words, $\mathbf{Y} = g(\mathbf{A})$ has the same distribution
    for each particular spherical distribution. An analogous situation is true
    for the class of vector-spherical distribution, see \citet[Section 5.1, pp. 154-156]{fz:90}
    and \citet[Section 5.1, pp. 182-189]{gv:93}.
    In particular note the following: let $\mathbf{A} \in \mathcal{L}_{m,n}^{\beta}$, be an
    $n \times m$ spherical random matrix $SS^{\beta}(n,m,h)$, such that
    $\mathbf{A} = (\mathbf{A}_{1}^{T}|\mathbf{A}_{2}^{T})^{T}$ with $\mathbf{A}_{i}
    \in \mathcal{L}_{m,n_{i}}^{\beta}$, $n_{i} \geq m$ $i = 1,2$ and $n = n_{1} +
    n_{2}$. Then
    \begin{enumerate}
      \item  the $T$ type II random matrix defined as
        $$
          \mathbf{T} = \mathbf{A}_{1}(\mathbf{A}_{2}^{*}\mathbf{A}_{2})^{-1/2},
        $$
      \item  the Gegenbauer type II random matrix defined as
        $$
          \mathbf{R} = \mathbf{A}_{1}(\mathbf{A}_{1}^{*}\mathbf{A}_{1} + \mathbf{A}_{2}^{*}
          \mathbf{A}_{2})^{-1/2},
        $$
      \item  the classical Jacobi random matrix defined as
        $$
          \mathbf{B} =(\mathbf{A}_{1}^{*}\mathbf{A}_{1} + \mathbf{A}_{2}^{*}\mathbf{A}_{2})^{-1/2}
          (\mathbf{A}_{1}^{*}\mathbf{A}_{1})(\mathbf{A}_{1}^{*}\mathbf{A}_{1} +
          \mathbf{A}_{2}^{*}\mathbf{A}_{2})^{-1/2},
        $$
      \item and the $F$ random matrix defined as
        $$
          \mathbf{F} = (\mathbf{A}_{2}^{*}\mathbf{A}_{2})^{-1/2}
          (\mathbf{A}_{1}^{*}\mathbf{A}_{1})(\mathbf{A}_{2}^{*}\mathbf{A}_{2})^{-1/2},
        $$
    \end{enumerate}
    are invariant under the class of spherical distributions, see \citet[Section 3.5, pp. 110-116]{fz:90}.
    Analogous results are concluded for the class of vector-spherical distributions,
    see \citet[Theorem 5.3.1, p.182]{gv:93}. Moreover, the classical Jacobi and $F$ random matrices have the
    same joint element distribution if $\mathbf{A}$ is a vector-spherical or spherical random
    matrix, see \citet[Theorems 3.5.1 and 3.5.5]{fz:90} and \citet[Theorem 5.3.1]{gv:93}.
  \item Observe that the Gegenbauer-Laguerre type II ensemble is indeed
    the classical Jacobi random matrix ensemble. Also note that, if $\mathbf{F}$ is
    an $F$ random matrix, then $\mathbf{B} = \mathbf{I}_{m} - (\mathbf{I}_{m} +
    \mathbf{F})^{-1}$, and if $\mathbf{B}$ is a classical Jacobi random matrix ensemble then
    $\mathbf{F} = (\mathbf{I}_{m} - \mathbf{B})^{-1} - \mathbf{I}_{m}$. Therefore $\mathbf{F}$ can
    be termed a \textbf{modified Jacobi random matrix ensemble}. Similarly, the $T$-Laguerre
    type II random matrix ensemble is indeed the modified Jacobi random matrix.
  \item There remain two final remarks:
    \begin{enumerate}
      \item As \citet{es:08} observed, the classical Jacobi random matrix can be
        obtained from the Fourier random matrix. Moreover, this procedure can be applied to any
        function $g(A)$, invariant under the class of spherical distributions.
        Form the statistician's point of view, this approach has been used by \citet{k:70} and
        \citet{c:96}, in the real case. The latter studied the $T$, inverted $T$, $F$
        and beta distributions.
      \item Finally, note that the matrix factorisation associated with classical and modified Jacobi ensembles
        could be the singular value decompositions, but applied to Gegenbauer and $T$ random
        matrices, respectively, see \citet{er:05}.
    \end{enumerate}
\end{enumerate}
\end{rem}

\subsection{Joint eigenvalue densities}\label{subsec34}

\begin{thm}\label{teovseigd}
Let $\mathbf{A}$ be an $m \times m$ matrix from the vector-spherical ensemble. Then
its joint eigenvalue density function is
\begin{equation}\label{vseigd}
    f^{\beta}_{\mathbf{\Lambda}}(\mathbf{\Lambda}) =
    \frac{c^{\beta}(m,n)2^{m(m-1)\beta/2}
    \pi^{\beta m^{2}/2 + \tau}}{\Gamma_{m}^{\beta}[\beta
    m/2]}\prod_{i <j}^{m}(\lambda_{i} - \lambda_{j})^{\beta} h\left(\sum_{i = 1}^{m}
    \lambda_{i}^{2}\right),
\end{equation}
where $\mathbf{\Lambda} = \diag(\lambda_{1}, \dots,\lambda_{m})$.
\end{thm}
\begin{proof}
From (\ref{vsed}) and Lemma \ref{lemsd}, the joint density of $\mathbf{W} \in
\mathfrak{U}^{\beta}(m)$ and $\mathbf{\Lambda} = \diag(\lambda_{1},
\dots,\lambda_{m})$, with $\mathbf{A} = \mathbf{W}\mathbf{\Lambda} \mathbf{W}^{*}$ is
$$
  2^{m(m-1)\beta/4}c^{\beta}(m,m)h(\tr \mathbf{\Lambda}^{2})2^{-m} \pi^{\tau}
   \prod_{i < j}^{m} (\lambda_{i} - \lambda_{j})^{\beta}
   (d\mathbf{\Lambda})(\mathbf{W}^{*}d\mathbf{W}).
$$
The marginal density desired is found by integrating over $\mathbf{W} \in
\mathfrak{U}^{\beta}(m)$ using (\ref{vol}).
\end{proof}

Similarly,
\begin{thm}\label{teosseigd}
Let $\mathbf{A}$ be an $m \times m$ matrix from the spherical ensemble. Then its
joint eigenvalue density function is
$$
    f^{\beta}_{\mathbf{\Lambda}}(\mathbf{\Lambda}) =
    \frac{c^{\beta}(m,n)2^{m(m-1)\beta/2}
    \pi^{\beta m^{2}/2 + \tau}}{\Gamma_{m}^{\beta}[\beta
    m/2]}\prod_{i <j}^{m}(\lambda_{i} - \lambda_{j})^{\beta}\qquad\qquad
$$\vspace{-.5cm}
\begin{equation}\label{sseigd}\hspace{7cm}
   h\left(\diag\left(\lambda_{1}^{2}, \dots,\lambda_{m}^{2}\right)\right)
\end{equation}
where $\mathbf{\Lambda} = \diag(\lambda_{1}, \dots,\lambda_{m})$.
\end{thm}
\begin{proof}
The proof is analogous to that given for (\ref{vseigd}); note however that
$h\left(\lambda\left(\mathbf{A}^{2}\right)\right) =
h\left(\lambda\left((\mathbf{W}\mathbf{\Lambda} \mathbf{W}^{*})^{2}\right)\right) =
h\left(\lambda\left(\mathbf{\Lambda}^{2}\right)\right) =
h\left(\diag(\lambda_{1}^{2}, \dots,\lambda_{m}^{2})\right)$.
\end{proof}

\begin{thm}\label{teovsleigd}
Let $\mathbf{S}$ be an $m \times m$ matrix from the vector-spherical-Laguerre
ensemble. Then its joint eigenvalue density function is
$$
  f^{\beta}_{\mathbf{\Lambda}}(\mathbf{\Lambda}) =
    \frac{c^{\beta}(m,n)\pi^{\beta m(n + m)/2 + \tau}}{\Gamma_{m}^{\beta}[\beta
    n/2]\Gamma_{m}^{\beta}[\beta m/2]} \prod_{i=1}^{m} \lambda^{\beta(n - m + 1)/2 - 1}
    \prod_{i <j}^{m}(\lambda_{i} - \lambda_{j})^{\beta}
$$\vspace{-.5cm}
\begin{equation}\label{vsleigd}\hspace{7cm}
     \times \ h\left(\sum_{i = 1}^{m} \lambda_{i}\right),
\end{equation}
where $\mathbf{\Lambda} = \diag(\lambda_{1}, \dots,\lambda_{m})$.
\end{thm}
\begin{proof}
Let $\mathbf{S} = \mathbf{W}\mathbf{\Lambda} \mathbf{W}^{*}$, from (\ref{vsld}) and
Lemma \ref{lemsd}, then the joint density of $\mathbf{W} \in \mathfrak{U}^{\beta}(m)$
and $\mathbf{\Lambda} = \diag(\lambda_{1}, \dots,\lambda_{m})$ is
$$
  \frac{c^{\beta}(m,n) \pi^{\beta mn/2}}{\Gamma_{m}^{\beta}[\beta n/2]}
  |\mathbf{S}|^{\beta(n - m + 1)/2 - 1} h(\tr(\mathbf{S}))2^{-m} \pi^{\tau}
  \prod_{i < j}^{m} (\lambda_{i} - \lambda_{j})^{\beta}
  (d\mathbf{\Lambda})(\mathbf{W}^{*}d\mathbf{W}).
$$
The marginal density desired is found by integrating over $\mathbf{W} \in
\mathfrak{U}^{\beta}(m)$ using (\ref{vol}).
\end{proof}

Similarly,
\begin{thm}\label{teossleigd}
Let $\mathbf{S}$ be an $m \times m$ matrix from the spherical-Laguerre ensemble. Then
its joint eigenvalue density function is
$$
  f^{\beta}_{\mathbf{\Lambda}}(\mathbf{\Lambda}) =
    \frac{c^{\beta}(m,n)\pi^{\beta m(n + m)/2 + \tau}}{\Gamma_{m}^{\beta}[\beta
    n/2]\Gamma_{m}^{\beta}[\beta m/2]} \prod_{i=1}^{m} \lambda^{\beta(n - m + 1)/2 - 1}
    \prod_{i <j}^{m}(\lambda_{i} - \lambda_{j})^{\beta}
$$\vspace{-.5cm}
\begin{equation}\label{ssleigd}\hspace{6cm}
    \times \  h\left(\diag\left(\lambda_{1}, \dots,\lambda_{m}\right)\right),
\end{equation}
where $\mathbf{\Lambda} = \diag(\lambda_{1}, \dots,\lambda_{m})$.
\end{thm}
Given the differential form of the joint element density function of the Fourier
ensemble, the corresponding joint eigenvalue density function cannot be obtained by
applying the Theorem \ref{teosseigd}. This density function has been obtained by
\citet[Lemma 10.4.4, p. 198]{me:91} and \citet[Proposition 2.3, p. 61]{f:05}. Based
on the wedge product, this is now obtained and
$\Vol\left(\mathfrak{P}^{\beta}_{S}(m)\right)$ is calculated indirectly, see
(\ref{vol2}).

Let $\mathbf{H} \in \mathfrak{P}^{\beta}_{S}(m)$, then there exist $\mathbf{U} \in
\mathfrak{P}^{\beta}(m)$ such that $\mathbf{H} = \mathbf{UEU}^{*}$, where $\mathbf{E}
= \diag(\exp(i\theta_{1}), \dots \exp(i\theta_{m}))$, $\exp(i\theta_{i})$ are $m$
complex numbers on the unit circle, see \citet{me:91} and \citet{f:05}. Then the
joint density function of $\theta_{1}, \dots, \theta_{m}$ is
$$
  c^{\beta}(m) \prod_{l \leq j} |\exp(i\theta_{l})-\exp(i\theta_{j})|^{\beta}, \quad
  \theta_{j} \in (-\pi,\pi).
$$
From (\ref{vol2})
$$
  dF_{\mathbf{H}}(\mathbf{H}) =
  \frac{1}{\Vol\left(\mathfrak{P}^{\beta}_{S}(m)\right)}(\mathbf{H}d\mathbf{H}).
$$
Let $\mathbf{H} = \mathbf{UEU}^{*}$, then by \citet[Exercise 2 (iii), p. 63]{f:05},
recalling that $\mathbf{H} = \mathbf{H}^{*}$
$$
  (\mathbf{H}d\mathbf{H})= \prod_{l \leq j} |\exp(i\theta_{l})-\exp(i\theta_{j})|^{\beta} \left ( \bigwedge_{l =
  1}^{m} \theta_{l} \right) (\mathbf{U}^{*}d\mathbf{U})
$$
Note that it is necessary to divide the measure by  $2^{-m} \pi^{\tau}$ to normalise
the arbitrary phases of the $m$ elements in the first row of $\mathbf{U}$, where
$\tau$ is given in Lemma \ref{lemsvd}. Then the joint density of $\diag(\theta_{1},
\dots \theta_{m})$ and $\mathbf{U}$ is
$$
  \frac{2^{-m} \pi^{\tau}}{\Vol\left(\mathfrak{P}^{\beta}_{S}(m)\right)}\prod_{l \leq j}
  |\exp(i\theta_{l})-\exp(i\theta_{j})|^{\beta} \left ( \bigwedge_{l = 1}^{m}
  \theta_{l} \right) (\mathbf{U}^{*}d\mathbf{U}).
$$
By integrating with respect to $\mathbf{U}$ we obtain the marginal of $\theta_{1},
\dots \theta_{m}$
\begin{equation}\label{fourier}
    \frac{2^{-m} \pi^{\tau} \Vol\left(\mathfrak{P}^{\beta}(m)\right)}{\Vol \left(
    \mathfrak{P}^{\beta}_{S}(m)\right)}\prod_{l \leq j}
    |\exp(i\theta_{l})-\exp(i\theta_{j})|^{\beta}  \bigwedge_{l = 1}^{m}
    \theta_{l}.
\end{equation}
From the Morris integral \citet[Equation (3.4), p. 122]{f:05}
\begin{eqnarray*}
  M_{m}[0,0, \beta/2] &=& \displaystyle(2\pi)^{-m} \int_{-\pi}^{\pi}\cdots \int_{-\pi}^{\pi} \prod_{l \leq j}
  |\exp(i\theta_{l})-\exp(i\theta_{j})|^{\beta} \bigwedge_{l = 1}^{m}  \theta_{l}\\
  &=& \displaystyle \frac{(2\pi)^{-m} \Gamma[\beta m/2 +1]}{(\Gamma[\beta/2 +
  1])^{m}},
\end{eqnarray*}
then
$$
  c^{\beta}(m) = \frac{\Gamma[\beta m/2 +1]}{(\Gamma[\beta/2 +
  1])^{m}}.
$$
However, from (\ref{fourier}) we have
$$
  c^{\beta}(m) = \frac{\Gamma[\beta m/2 +1]}{(\Gamma[\beta/2 +
  1])^{m}} = \frac{2^{-m} \pi^{\tau} \Vol\left(\mathfrak{P}^{\beta}(m)\right)}{\Vol \left(
  \mathfrak{P}^{\beta}_{S}(m)\right)},
$$
Therefore
$$
  \Vol \left(\mathfrak{P}^{\beta}_{S}(m)\right) = \frac{2^{m} \pi^{\beta m^{2}/2 + \tau}
  (\Gamma[\beta/2 +1])^{m}}{\Gamma_{m}^{\beta}[\beta m/2] \Gamma[\beta m/2 +1]}.
$$

Table \ref{table2} summarises the kernels of the joint eigenvalue density function of
all the particular ensembles studied in the subsection \ref{subsec33}, using the
names proposed in Remark \ref{rem33}. The corresponding normalisation constants are
obtained from their joint element density functions and from Theorems
\ref{teovseigd}, \ref{teosseigd}, \ref{teovsleigd} and \ref{teossleigd}.

\begin{table}
  \centering
  \caption{Kernels of the joint eigenvalue density function of all the particular ensembles}\label{table2}
  \begin{tabular}{l|l}
    \hline \hline
    $\build{\mbox{Ensemble}}{}{\lambda_{i} \in (a,b)}$ & Kernel \\
    \hline\hline
    $\build{\mbox{Hermite}}{}{\lambda_{i} \in (-\infty, \infty)}$
    & $\displaystyle \exp\left\{\frac{\beta}{2} \sum_{i = 1}^{m}\lambda_{i}^{2}\right\}
    \prod_{i <j}^{m}(\lambda_{i} - \lambda_{j})^{\beta}$ \\
    $\build{T\mbox{ type I}}{}{\lambda_{i} \in (-\infty, \infty)}$
    & $\displaystyle \left(1 +  \sum_{i = 1}^{m}\lambda_{i}^{2}\right)
    ^{-\beta[m(m-1)/2 + \nu]/2 - m/2} \prod_{i <j}^{m}(\lambda_{i} - \lambda_{j})^{\beta}$ \\
    $\build{\mbox{Gegenbauer type I}}{}{\lambda_{i} \in (-1, 1)}$ &
    $\displaystyle  \left(1 - \sum_{i = 1}^{m}\lambda_{i}^{2}\right)^{\beta q}
    \prod_{i <j}^{m}(\lambda_{i} - \lambda_{j})^{\beta}$ \\
    $\build{T\mbox{ type II}}{}{\lambda_{i} \in (-\infty, \infty)}$
    & $\displaystyle \prod_{i = 1}^{m} (1 + \lambda_{i}^{2})^{-\beta(n + \nu)/2}
    \prod_{i <j}^{m}(\lambda_{i} - \lambda_{j})^{\beta}$ \\
    $\build{\mbox{Gegenbauer type II}}{}{\lambda_{i} \in (-1, 1)}$
    & $\displaystyle \prod_{i = 1}^{m} (1-\lambda_{i}^{2})^{\beta(\nu - m + 1)/2 - 1}
    \prod_{i <j}^{m}(\lambda_{i} - \lambda_{j})^{\beta}$ \\
    $\build{\mbox{Laguerre}}{}{\lambda_{i} \in (0, \infty)}$
    & $\displaystyle \prod_{i = 1}^{m}\lambda_{i}^{\beta(n - m + 1)/2 - 1}
    \exp\left\{\frac{\beta}{2} \sum_{i = 1}^{m}\lambda_{i}\right\}
    \prod_{i <j}^{m}(\lambda_{i} - \lambda_{j})^{\beta}$ \\
    $\build{\mbox{Modified Jacobi type I}}{}{\lambda_{i} \in (0, \infty)}$
    & $\displaystyle \prod_{i = 1}^{m}\lambda_{i}^{\beta(n - m + 1)/2 - 1}
    \left(1+ \sum_{i = 1}^{m}\lambda_{i} \right)^{-\beta(n + \nu)/2}
    \prod_{i <j}^{m}(\lambda_{i} - \lambda_{j})^{\beta}$ \\
    $\build{\mbox{Jacobi type I}}{}{\lambda_{i} \in (0, 1)}$
    & $\displaystyle \prod_{i = 1}^{m}\lambda_{i}^{\beta(n - m + 1)/2 - 1} \left(1-
    \sum_{i = 1}^{m}\lambda_{i}\right)^{\beta(\nu - m + 1)/2 - 1}
    \prod_{i <j}^{m}(\lambda_{i} - \lambda_{j})^{\beta}$ \\
    $\build{\mbox{Modified Jacobi}}{}{\lambda_{i} \in (0, \infty)}$
    & $\displaystyle \prod_{i = 1}^{m}\lambda_{i}^{\beta(n - m + 1)/2 - 1}
    (1+ \lambda_{i})^{-\beta(n + \nu)/2} \prod_{i <j}^{m}(\lambda_{i} - \lambda_{j})^{\beta}$ \\
    $\build{\mbox{Jacobi}}{}{\lambda_{i} \in (0, 1)}$
    & $\displaystyle \prod_{i = 1}^{m}\lambda_{i}^{\beta(n - m + 1)/2 - 1} (1-
    \lambda_{i})^{\beta(\nu - m + 1)/2 - 1} \prod_{i <j}^{m}(\lambda_{i} - \lambda_{j})^{\beta}$ \\
    $\build{\mbox{Fourier}}{}{\lambda_{i} \in (-\pi, \pi)}$
    & $\displaystyle \prod_{l \leq j}^{m}|\exp(i\theta_{l}) - \exp(i\theta_{j})|^{\beta}$ \\
    \hline\hline
  \end{tabular}
\end{table}

Finally, another interesting property of the vector-spherical random matrix ensemble
is that its joint normalised eigenvalue density function is invariant under all
classes of vector-spherical distribution, where the normalised eigenvalue can be
defined as $\delta_{i} = \lambda_{i}/r$. In particular, $r$ can be defined as
$$
  r = \left(\sum_{i = 1}^{m} \lambda^{2} \right)^{1/2}.
$$
This problem has been studied in the context of shape theory and the corresponding
density function was obtained by \citet{dggjrq:03}.

\section{Multivariate statistics}\label{sec4}
Many areas of multivariate statistics (and of statistics in general), such as
multiple time series or econometrics, can be enriched with the use of the tools and
ideas of random matrix theory, see \citet{ha:70}. However, in this section these
tools and ideas are used to develop a unified theory of multivariate distributions
for normed division algebras.

\begin{thm}\label{teoled}
Let $\mathbf{Z} \in \mathcal{L}_{m,n}^{\beta}$ be a matrix variate left elliptical
distribution with density function
$$
  c^{\beta}(m,n)h(\mathbf{Z}^{*}\mathbf{Z}).
$$
Therefore if $\mathbf{X} = \mathbf{AZB} + \boldsymbol{\mu}$, with $\mathbf{A} \in
\mathcal{L}_{n,n}^{\beta}$, $\mathbf{B} \in \mathcal{L}_{m,m}^{\beta}$ and
$\boldsymbol{\mu} \in \mathcal{L}_{m,n}^{\beta}$, constant matrices such that
$\mathbf{A}^{*}\mathbf{A} = \mathbf{\Theta} \in \mathfrak{P}_{n}^{\beta}$ and
$\mathbf{B}^{*}\mathbf{B} = \mathbf{\Sigma} = (\mathbf{\Sigma}^{1/2})^{2} \in
\mathfrak{P}_{m}^{\beta}$, then
$$
  \frac{c^{\beta}(m,n)}{|\mathbf{\Sigma}|^{\beta n/2}|\mathbf{\Theta}|^{\beta m/2}}h\left (\mathbf{\Sigma}^{-1/2}
  (\mathbf{X} - \boldsymbol{\mu})^{*} \mathbf{\Theta}^{-1}(\mathbf{X} - \boldsymbol{\mu})
  \mathbf{\Sigma}^{-1/2}\right ).
$$
where $c^{\beta}(m,n)$ is a normalisation constant.
\end{thm}
\begin{proof}
The proof is analogous to that given in the real case, but considering the Jacobian
of the transformation $\mathbf{Y} = \mathbf{AXB} + \boldsymbol{\mu}$ defined by Lemma
\ref{lemlt}.
\end{proof}
\begin{cor}
Under the hypothesis of Theorem \ref{teoled}:
\begin{enumerate}
  \item for the matrix variate spherical elliptical distribution, its density function
  is
  $$
    \frac{c^{\beta}(m,n)}{|\mathbf{\Sigma}|^{\beta n/2}|\mathbf{\Theta}|^{\beta m/2}}
    h\left ( \lambda \left (\mathbf{\Sigma}^{-1} (\mathbf{X} - \boldsymbol{\mu})^{*}
    \mathbf{\Theta}^{-1}(\mathbf{X} - \boldsymbol{\mu})\right )\right ),
  $$
  this fact being denoted as $\mathbf{X} \sim \mathcal{SSE}_{n \times m}^{\beta}(\boldsymbol{\mu},
  \mathbf{\Sigma}, \mathbf{\Theta}, h)$.
  \item for the matrix variate vector-spherical elliptical distribution, its density function
  is
  $$
    \frac{c^{\beta}(m,n)}{|\mathbf{\Sigma}|^{\beta n/2}|\mathbf{\Theta}|^{\beta m/2}}h\left
    ( \tr\left (\mathbf{\Sigma}^{-1}(\mathbf{X} - \boldsymbol{\mu})^{*} \mathbf{\Theta}^{-1}
    (\mathbf{X} - \boldsymbol{\mu})\right )\right ),
  $$
  this fact being denoted as $\mathbf{X} \sim \mathcal{VSE}_{n \times m}^{\beta}(\boldsymbol{\mu},
  \mathbf{\Sigma}, \mathbf{\Theta}, h)$.
  \item and for the matrix variate normal distribution, its density function
  is (taking $h(u) = \exp(-\beta u/2)$)
  $$
  \frac{1}{\left(2\pi \beta^{-1}\right)^{\beta mn/2}|\mathbf{\Sigma}|^{\beta n/2}
  |\mathbf{\Theta}|^{\beta m/2}}\etr\left \{- \frac{\beta}{2} \mathbf{\Sigma}^{-1}
  (\mathbf{X} - \boldsymbol{\mu})^{*} \mathbf{\Theta}^{-1}(\mathbf{X} - \boldsymbol{\mu})
  \right \},
  $$
  denoting this fact as $\mathbf{X} \sim \mathcal{N}_{n \times m}^{\beta}(\boldsymbol{\mu},
  \mathbf{\Sigma}, \mathbf{\Theta})$.
\end{enumerate}
\end{cor}
Many other particular vector-spherical or spherical elliptical distributions can be
obtained by simply specifying the function $h(\cdot)$ in a similar way to Corollaries
\ref{corvspd} and \ref{corsspd}.

Elliptical and, in particular, normal symmetric random matrices have received less
attention in multivariate statistics, but analogous results can be obtained in a
similar way using Lemma \ref{lemhlt}.
\begin{thm}\label{teogwd}
Let us define $\mathbf{S} = \mathbf{X}^{*} \mathbf{\Theta}^{-1}\mathbf{X} \in
\mathfrak{P}_{m}^{\beta}$, then
\begin{enumerate}
  \item if $\mathbf{X} \sim \mathcal{SSE}_{n \times m}^{\beta}(\boldsymbol{0},
  \mathbf{\Sigma}, \mathbf{\Theta}, h)$, the density function of $\mathbf{S}$ is
  $$
    \frac{d^{\beta}(m,n) \pi^{\beta mn/2}}{\Gamma_{m}^{\beta}[\beta n/2]|\mathbf{\Sigma}|^{\beta n/2}}
    |\mathbf{S}|^{\beta(n - m + 1)/2 - 1} h(\lambda(\mathbf{\Sigma}^{-1}\mathbf{S})).
  $$
  This distribution is known as the spherical-generalised-Wishart distribution and it is
  denoted as $\mathbf{S} \sim \mathcal{SGW}_{m}^{\beta}(n, \mathbf{\Sigma}, h)$.
  \item if $\mathbf{X} \sim \mathcal{VSE}_{n \times m}^{\beta}(\boldsymbol{0},
  \mathbf{\Sigma}, \mathbf{\Theta}, h)$, the density function of $\mathbf{S}$ is
  $$
    \frac{d^{\beta}(m,n) \pi^{\beta mn/2}}{\Gamma_{m}^{\beta}[\beta n/2]|\mathbf{\Sigma}|^{\beta n/2}}
    |\mathbf{S}|^{\beta(n - m + 1)/2 - 1} h(\tr(\mathbf{\Sigma}^{-1}\mathbf{S})).
  $$
  This distribution is known as the vector-spherical-generalised-Wishart distribution and it is
  denoted as $\mathbf{S} \sim \mathcal{VSGW}_{m}^{\beta}(n, \mathbf{\Sigma}, h)$.
\end{enumerate}
\end{thm}
\begin{proof}
Let $\mathbf{S} = \mathbf{X}^{*} \mathbf{\Theta}^{-1}\mathbf{X}$, and note that
$$
  \mathbf{Y} = \mathbf{\Theta}^{-1/2}\mathbf{X} \sim \mathcal{SSE}_{n \times
    m}^{\beta}(\boldsymbol{0}, \mathbf{\Sigma}, \mathbf{I}_{n}, h)
$$
with $\left(\mathbf{\Theta}^{-1/2}\right)^{2} = \mathbf{\Theta}$. Therefore
$\mathbf{S} = \mathbf{X}^{*} \mathbf{\Theta}^{-1}\mathbf{X} =
\mathbf{Y}^{*}\mathbf{Y}$. The desired results are follow from (\ref{sld}) and
(\ref{vsld}), respectively.
\end{proof}

\begin{cor}
Under the hypothesis of Theorem \ref{teogwd} it follows that,
\begin{enumerate}
  \item the Wishart random matrix has a density function
    $$
      \frac{1}{\left(2\beta^{-1}\right)^{\beta mn/2} \Gamma_{m}^{\beta}[\beta mn/2] |\mathbf{\Sigma}|^{\beta n/2}}
      |\mathbf{S}|^{\beta(n - m +1)/2 -1}\etr\{-\beta \mathbf{\Sigma}^{-1}\mathbf{S}/2\}
    $$
    where $n \geq (m-1)\beta$ and this is denoted as $\mathbf{S} \sim \mathcal{W}_{m}^{\beta}(n, \mathbf{\Sigma})$.;
  \item the matrix variate beta type I distribution has a density function
    $$
      \frac{1}{\mathcal{B}_{m}^{\beta}[\beta n/2,\beta \nu/2]} |\mathbf{S}|^{\beta(n - m +1)/2 -1}
      |\mathbf{I}_{m} - \mathbf{S}|^{\beta(\nu - m +1)/2 -1},
    $$
    where $\nu  \geq (m-1)\beta$ and $n \geq (m-1)\beta$ and this is
  denoted as $\mathbf{S} \sim \mathcal{BI}_{m}^{\beta}(n, \nu)$;
  \item and the matrix variate beta type II distribution has a density function
    $$
      \frac{1}{\mathcal{B}_{m}^{\beta}[\beta m/2,\beta n/2]} |\mathbf{S}|^{\beta(n - m +1)/2 -1}
      |\mathbf{I}_{m} + \mathbf{S}|^{-\beta(n + \nu)/2},
    $$
    where $\nu \geq (m-1)\beta$ and $n \geq (m-1)\beta$ and this is
  denoted as $\mathbf{S} \sim \mathcal{BII}_{m}^{\beta}(n, \nu)$.
\end{enumerate}
\end{cor}
\begin{proof}
The desired results are obtained from Theorem \ref{teogwd} directly, after defining
$h(\tr(\mathbf{\Sigma}^{-1}\mathbf{S})) = \etr\{-\beta
\mathbf{\Sigma}^{-1}\mathbf{S}/2\}$, $h(\lambda(\mathbf{\Sigma}^{-1} \mathbf{S})) =
|\mathbf{I}_{m} - \mathbf{S}|^{\beta(\nu - m +1)/2 -1}$ ($\mathbf{\Sigma} =
\mathbf{I}_{m}$) and $h(\lambda(\mathbf{\Sigma}^{-1} \mathbf{S})) = |\mathbf{I}_{m} +
\mathbf{S}|^{-\beta(n + \nu)/2}$ ($\mathbf{\Sigma} = \mathbf{I}_{m}$), respectively.
\end{proof}

\section*{Conclusions}

Although most results about Jacobians are known in the context of random matrix
theory, in general for statisticians they are less familiar, as are the corresponding
technical tools in the context of normed division algebras. Thus the importance of
addressing this area of study.

In the field of random matrix theory, various classes of ensembles are proposed,
including many of those studied previously, see \citet[Section 4.1.4, p. 177]{f:05}.
What is most important is that these new classes of ensembles contain many other
ensembles of potential interest, which may enable us to study phenomena and
experiments under more general conditions.

Analogously to current random matrix theory, see \citet{er:05}, in the present study
we propose a unified theory of matrix variate distribution for normed division
algebra, that is, for real, complex, quaternion, and octonion cases.

\section*{Acknowledgements}
This research work was partially supported  by CONACYT-M\'exico, Research Grant No. \
81512 and IDI-Spain, Grants No. FQM2006-2271 and MTM2008-05785. This paper was
written during J. A. D\'{\i}az- Garc\'{\i}a's stay as a visiting professor at the
Department of Statistics and O. R. of the University of Granada, Spain.

\end{document}